\PassOptionsToPackage{unicode}{hyperref}
\PassOptionsToPackage{hyphens}{url}
\documentclass[
]{article}
\usepackage{amsmath,amssymb}
\usepackage{lmodern}
\usepackage{iftex}
\usepackage{varwidth}
\usepackage[pdftex,linktocpage]{hyperref}
\usepackage{amsthm}
\ifPDFTeX
  \usepackage[T1]{fontenc}
  \usepackage[utf8]{inputenc}
  \usepackage{textcomp} 
\else 
  \usepackage{unicode-math}
  \defaultfontfeatures{Scale=MatchLowercase}
  \defaultfontfeatures[\rmfamily]{Ligatures=TeX,Scale=1}
\fi
\IfFileExists{upquote.sty}{\usepackage{upquote}}{}
\IfFileExists{microtype.sty}{
  \usepackage[]{microtype}
  \UseMicrotypeSet[protrusion]{basicmath} 
}{}
\makeatletter
\@ifundefined{KOMAClassName}{
  \IfFileExists{parskip.sty}{%
    \usepackage{parskip}
  }{
    \setlength{\parindent}{0pt}
    \setlength{\parskip}{6pt plus 2pt minus 1pt}}
}{
  \KOMAoptions{parskip=half}}
\makeatother
\usepackage{xcolor}
\usepackage[normalem]{ulem}
\ifLuaTeX
  \usepackage{selnolig}  
\fi
\IfFileExists{bookmark.sty}{\usepackage{bookmark}}{\usepackage{hyperref}}
\IfFileExists{xurl.sty}{\usepackage{xurl}}{} 
\hypersetup{
    colorlinks=true,
    linkcolor=blue,
    filecolor=magenta,      
    urlcolor=cyan,
    pdftitle={Overleaf Example},
    pdfpagemode=FullScreen,
    }

\author{}
\date{}

\begin{document}

\begin{center}
    
\section{{An introduction to univalent function theory and the Bieberbach conjecture}}

\subsection{Jiakai Qu}\footnote{25 Gordon Street (UCL Union Building), zcahqux@ucl.ac.uk} \footnote{\emph{2020 MSC 30C50, 30C55}}

~

University College London

December 22, 2022

\end{center}

~ ~ ~ ~ ~

\textbf{Abstract}

The purpose of this paper is to make an introduction to univalent
function theory for readers of any level, assuming only foundational knowledge in real and complex analysis. In particular, we state and proof (with details) important theorems utilized in proving Bieberbach's conjecture, especially those that were missed or merely sketched by other texts on univalent
functions. We will finally prove the conjecture following the proof given by Lenard Weinstein in a comprehensive and self-contained manner, exposing its full technical details in each step. Readers are encouraged to make use of this paper to accompany the preliminary chapters of Duren's \emph{Univalent Functions} and Lenard Weinstein's \emph{the Bieberbach Conjecture}. 

~

\textbf{Keywords:} Univalent functions, Classical Loewner theory,
Bieberbachs's conjecture

~

\textbf{Acknowledgement:} The author would like to thank Professor Alexander Sobolev (University College
London), for his patience, kindness, and exceptional expertise,
without which this paper could not have possibly been finished.

\newpage


\section {Contents}

\hypertarget{Introduction}{%
\section{0 Introduction}\label{Introduction}}

\hypertarget{univalent-functions}{%
\section{1 Univalent Functions}\label{univalent-functions}}

\begin{quote}
1.1 Prerequisites

1.2 Univalent functions

1.3 Preliminary Results

1.4 Further Theorems
\end{quote}

\hypertarget{the-bieberbach-conjecture}{%
\section{2 The Bieberbach conjecture}\label{the-bieberbach-conjecture}}

\begin{quote}
2.1 The Bieberbach conjecture

2.2 Robertson's conjecture

2.3 Milin's conjecture
\end{quote}

\hypertarget{lowners-theory}{%
\section{3 Lowner's theory}\label{lowners-theory}}

\begin{quote}
3.1 Caratheodory's kernel convergence theorem

3.2 Lowner's theory and Further lemmas
\end{quote}

\hypertarget{proof-of-bieberbachs-conjecture}{%
\section{4 Proof of Bieberbach's
conjecture}\label{proof-of-bieberbachs-conjecture}}

\begin{quote}
4.1 Legendre Polynomials

4.2 Remarks on proven results

4.2 Weinstein's Proof of the conjecture
\end{quote}

\newpage

~

\begin{center}
    
~\\

~\\

\emph{To Bangwen and Huimin}

\end{center}

\newpage

\textbf{
\section{0 Introduction}
}

Consider in the unit disc\(\mathbb{\ D\ }\)the Taylor series expansion of
univalent functions\(\ f\), normalized such that\(\ f(0) = 0\)
and\(\ f^{'}(0) = 1\):

\[f(z) = z + \sum_{n = 2}^{\infty}{a_{n}z^{n}} \]

This is the standard definition of the family\(\ S\ \)of univalent functions (see \hyperlink{D1}{[D1]}, p.26 or \hyperlink{Def1.2.1}{Def 1.2.1} of this paper); by the Riemann Mapping theorem the behaviours of this family extends to any univalent functions defined on arbitrary domains. In 1916, Ludwig Bieberbach proposed the later well-known \emph{Bieberbach's conjecture}

\fbox{\parbox{1.0\textwidth}{\ 
\emph{The coefficients of each function}\(\ f \in S\ \)\emph{satisfies}\(\ \left| a_{n} \right| \leq n\ \ \)\emph{for}\(\ \ n = 2,3,\ldots\)}}

Later Littlewood proved that (see \hyperlink{D1}{[D1]}, p.37 or \hyperlink{Thm2.1.2}{Thm 2.1.2} of this paper) by considering \(M_{p}(r,f) = \left\{ \frac{1}{2\pi}\int_{0}^{2\pi}\left| f\left( re^{i\theta} \right) \right|^{p}d\theta \right\}^{\frac{1}{p}}\ \)and
the growth theorem  we have that \(\ \left| a_{n} \right| \leq e \cdot n\). In 1936, M.S.Robertson
conjectured that (see \hyperlink{Rob}{Thm 2.2.1}) the coefficients of the function

\[g(z) = \sqrt{f\left( z^{2} \right)} = z + c_{3}z^{3} + c_{5}z^{5} + \ldots\]

satisfy

\[\sum_{k = 1}^{n}{\left| c_{2k - 1} \right|^{2} \leq n}\]

and showed that this is equivalent to Bieberbach's conjecture. If we consider 
\( g(z) = \sqrt{f\left( z^{2} \right)} \) and define the logarithmic coefficients \( \gamma_{n} \) by

\[\log\frac{g\left( \sqrt{z} \right)}{\sqrt{z}} = \frac{1}{2}\log\left( \frac{f(z)}{z} \right) = \sum_{n = 1}^{\infty}{\gamma_{n}z^{n}}\]

then it can be shown that the following conjectured inequality is also equivalent to Bieberbach's conjecture, as proposed by I.M.Milin in 1964

\[\sum_{m = 1}^{n}{\sum_{k = 1}^{m}\left( k\left| \gamma_{k} \right|^{2} - \frac{1}{k} \right)} \leq 0,\ n = 1,2,3\ldots\]

based on the well-known Lebedev-Milin Inequality (see \hyperlink{LM}{Thm 2.3.1}). 

After de Branges' first successful proof (see \hyperlink{Br1}{[Br1]}) of
the Bieberbach conjecture in 1985, Lenard Weinstein presented in 1991
a much shorter proof utilizing Legendre Polynomials and classical
Loewner theory, and this proof was famous for its conciseness, having
only 4 pages in length. 

There has been studies linking the two respective proofs by de Branges and Weinstein, for instance \hyperlink{Wi1}{[Wi1]}, as well as established survey papers on the comprehensive history and developments of Bieberbach's conjecture, such as \hyperlink{KO1}{[KO1]}.

In this paper we examine Weinstein's proof and results leading to it. Particularly, in Chapters 1 and 2, we first study the basic theory of univalent functions following the footpath of \hyperlink{D1}{[D1]} while supplementing proves to overlooked assumptions. In Chapter 3, we investigate Caratheodory's kernel convergence theorem and use it to derive Loewner's
PDE

\[p(z,s,t) = \frac{\frac{\partial f_{s}(z)}{\partial s}}{z\frac{\partial f_{s}(z)}{\partial z}}\]

where the function

\[f_{t}(z) = e^{t}\left\{ z + \sum_{n = 2}^{\infty}{b_{n}(t)z^{n}} \right\}\]

parameterizes single-slit mappings that converges uniformly to\(\ f \in S\).

In Chapter 4 we reveal the full details of Weinstein's proof of the
Bieberbach conjecture, after presenting a proof of the addition theorem of Legendre
Polynomials (which plays a central role in Weinstein's proof) in the form

\[P_{n}\left( \cos\theta_{1} \right)P_{n}\left( \cos\theta_{2} \right) + 2\sum_{k = 1}^{n}{( - 1)^{k}P_{n}^{- k}\left( \cos\theta_{1} \right)P_{n}^{k}\left( \cos\theta_{2} \right)}\cos(k\phi)\]

\[=P_{n}\left( \cos{\theta_{1}\cos{\theta_{2} + \sin{\theta_{1}\sin{\theta_{2}\cos\phi}}}} \right)\]

Then, after thoroughly studying Weinstein's computations, examine the generating function of Legendre Polynomials and uncover its link with the Koebe function, thus proving Milin's conjecture and henceforth Bieberbach's conjecture. In particular, we focus on unveiling the full computational details in \hyperlink{W1}{[W1]}: it has been commented by some that once exposed in full length, Weinstein's proof is much longer than expected; for instance J.H.Conway in \emph{Functions of One Complex Variable II} remarked: "Another proof can be found in Weinstein [1991].  Though this is shorter in the published form ..., it actually becomes longer if the same level of detail is provided." 

\newpage

\textbf{
\section{1 Univalent Functions}
}

\subsection{1.1 Prerequisites and notations}

The study of univalent function theory is a branch of complex analysis.
More precisely, it is part of geometric function theory. The basic idea
is to study the link between graphical transformations and one-to-one
complex-valued functions (in particular we only study functions with a
single variable in this text). To do so, we first make some remarks on
fundamental results in complex analysis, and introduce some theorems
useful for further studies. Readers are assumed to be familiar with: the
definition of holomorphic functions, conformal mappings, complex
sequences and series, and preliminary results in complex analysis, such
as Cauchy's theorem, Cauchy's Integral Formula, and Taylor/Laurent expansions etc. Everything begins from the basic definition of
univalence:

\uline{Definition 1.1.1}

\fbox{\parbox{1.0\textwidth}{\emph{
Let\(\ f\ \)be holomorphic on a domain\(\ D.\ \)Then we
call\(\ f\ \)\emph{univalent},or \emph{simple} if it is
one-to-one (injective) on\(\ D.\)
} }}

~

The following theorems are well-known and will not be proved in this
paper. Materials containing detailed proves of each theorems are
specified. Readers who are already familiar with the following results
may skip to the next section.

~

\hypertarget{Hurwitz}{\uline{Theorem 1.1.2 (Hurwitz's theorem)}}

\fbox{\parbox{1.0\textwidth}{\ 
\emph{Let} \(f_{n}\ \)\emph{be a sequence of analytic functions, not
equivalent to zero, in a domain}\(\ D.\ \)
~
\emph{Then if}\(\ f_{n} \rightarrow f\ \)\emph{uniformly on every compact
subset of}\(\ D,\ \)\emph{every zero of the limit
function}\(\ f\ \)\emph{is an accumulation point of a sequence of zeroes
of the functions}\(\ f_{n}\)\emph{.}
}}

For a proof, consult \hyperlink{C1}{[C1]} p.173 Theorem 7.5, or
\hyperlink{D1}{[D1]} p.4

\uline{Theorem 1.1.3}

\fbox{\parbox{1.0\textwidth}{\ 
\emph{Let} \(f_{n}\ \)\emph{be a sequence of analytic and univalent
functions in a domain}\(\ D.\ \)

\emph{Suppose that} \(f_{n} \rightarrow f\ \)\emph{uniformly on each
compact subset of}\(\ D.\ \)\emph{Then f is either univalent or constant
in}\(\ D.\ \)
}}

For a proof, consult \hyperlink{C1}{[C1]} p.174 Theorem 7.9, or
\hyperlink{D1}{[D1]} p.5
\newpage
\uline{Definition 1.1.4 (Normal Family)}

\fbox{\parbox{1.0\textwidth}{\ 
A family of analytic functions\(\mathcal{\ F\ }\)is called a
\emph{normal family} if every sequence of functions\(\ f_{n}\)
in\(\mathcal{\ F\ }\)contains a subsequence which converges uniformly on
every compact subset of\(\ D.\ \)
}}

The above definition is well-known and can be found in most texts on
complex analysis. In particular, readers are encouraged to consult
\hyperlink{S2}{[S2]} p.225 Ch.3.2, for detailed discussion on the concept of normal
families, complete exhaustion, and the proof of the following two
important theorems:

\hypertarget{Montel}{\uline {Theorem 1.1.5 (Montel's theorem)}}

\fbox{\parbox{1.0\textwidth}{\ 
\emph{If a family of analytic functions is locally bounded, then it is
normal.}

~

\emph{That is, if a family} \(\mathcal{F \in}D\ \)\emph{of analytic
functions satisfies} \(|f| \leq M\mathbb{\in R\ \ }\)\emph{for all
functions}\(\ \)\(f\mathcal{\in F,\ }\)\emph{then every sequence of functions}\(\ f_{n}\) \emph{in}\(\mathcal{\ \ F\ }\)\emph{contains a
sub-sequence that converges uniformly on every compact subset
of}\(\ D.\ \)

}}

\emph{The converse is also true: every normal family is locally
bounded.}

\uline {Theorem 1.1.6 (The Riemann mapping theorem)}

\fbox{\parbox{1.0\textwidth}{\ 
\emph{For any simply connected
domain}\(\ D\mathbb{\neq C\ }\)\emph{and}\(\ z_{0} \in D,\ \)\emph{there
exists a unique conformal mapping}\(\ \varphi\ \)\emph{from}\(\ D\ \)\emph{onto the unit disc, with}\(\ \varphi\left( z_{0} \right) = 0\ \)\emph{and}\(\ \varphi^{'}\left( z_{0} \right) > 0\)\emph{.}}}

There are other formulations for the proof of the above theorem. Other
than \hyperlink{S2}{[S2]}, readers are also encouraged to consult other approaches,
for example \hyperlink{D1}{[D1]} p.11, \hyperlink{B1}{[B1]} Ch.14 etc.

\uline{Theorem 1.1.7 (Vitali's theorem)}

\fbox{\parbox{1.0\textwidth}{\ 
\emph{Let the functions}\(\ f_{n}\ \)\emph{be analytic and locally
bounded in a domain}\(\ D,\ \)\emph{and suppose that}
\(\left\{ f_{n} \right\}\ \)\emph{converges pointwise on a set which has
an accumulation point in}\(\ D.\ \)

\emph{Then} \(\left\{ f_{n} \right\}\ \)\emph{converges uniformly on
every compact subset of}\(\ D.\ \)
}}

For a proof and further discussions on Montel's and Vitali's theorem,
readers are encouraged to consult \hyperlink{S1}{[S1]} p.227 Ch.6 or \hyperlink{C1}{[C1]} p.187
Theorem 7.36, or {[}D1{]} p.7-9.
\subsection{1.2 Univalent functions}

We mostly follow the footpath of {[}D1{]} for this section. Readers are
strongly encouraged to read this chapter alongside of {[}D1{]} as a
supplementary material.

The study of univalent functions concerns mostly of a particular
class\(\ S\), which simplifies certain problems by giving a straight
forward normalization:
\newpage
\hypertarget{Def1.2.1}{\uline{Definition 1.2.1}}

\fbox{\parbox{1.0\textwidth}{\ 
The class\(\ S\ \)is a family of functions\(\ f\ \), analytic and
univalent in the unit disk\(\mathbb{\ D,\ }\)normalized by the conditions\(\ f(0) = 0\ \)and\(\ f^{'}(0) = 1.\ \)}}

In view of the Riemann mapping theorem, when we wish to study some
arbitrary univalent function\(\ g\ \)that maps an arbitrary simply
connected domain\(\ D\ \)onto some domain\(\ U\ \)containing zero,
and\(\ g\left( z_{0} \right) = 0,\ \)we can always find an unique
conformal map\(\ \varphi:D \rightarrow \mathbb{D}\) and an unique
univalent map\(\ f:\mathbb{D \rightarrow}U\) in the family \(S\ \)such
that

\[g = f \circ \varphi\]

since we can fix the point\(\ z_{0} \in D\ \)and
let\(\ \varphi\left( z_{0} \right) = f(0) = 0.\ \)If\(\ U\ \)does not
contain the origin, then we can apply a transformation to reach the same
result. For any univalent \(p(z)\) within the unit
disc\(\mathbb{\ D\ }\)it holds that

\[f(z) = \frac{p(z) - p(0)}{p^{'}(0)}\]

is a member of \(S.\ \)Indeed, \(f(0) = 0\) and
\(f^{'}(0) = \frac{p^{'}(0)}{p^{'}(0)} = 1\).

For this reason, it is enough for us to study the class\(\ S\ \)and
translate the results to more generalized geometric theorems.

By definition of the class\(\ S\ \)all\(\ f \in S\ \)has Taylor series
expansion of the form

\[f(z) = z + a_{2}z^{2} + \ldots,\ |z| < 1\]

For example,

\[k(z) = z + 2z^{2} + 3z^{3} + \ldots\]

is a member of\(\ S.\ \)This function is called the Koebe function and
will be frequently referred to in future studies. Note that we can write

\[k(z) = \frac{z}{(1 - z)^{2}}\]

by manipulation on the geometric series \(1/(1 - z)\ \)(differentiate
and time by\(\ z\)).

We can also write this in the form

\[k(z) = \frac{z}{(1 - z)^{2}} = \frac{1}{4}\left( \frac{1 + z}{1 - z} \right)^{2} - \frac{1}{4}\]

which will again come to use in future arguments. Following the footpath of {[}D1{]}, we verify the following results
\hypertarget{Ele}{(\emph{\uline{Elementary transformations)}}}:

The class \(S\) is preserved under the following transformations if
\(f(z) \in S\)

(i) Conjugation: \(\overline{f\left( \overline{z} \right)} \in S\)

\begin{proof}
    
Notice that

\[\overline{f\left( \overline{z} \right)} = z + \overline{a_{2}}z^{2} + \overline{a_{3}}z^{3} + \ldots\]

Indeed, \(\overline{f(0)} = 0,\ \)and
\(\overline{f^{'}(0)} = 1 + 0 + 0\ldots = 1\)

Univalence:

\[\overline{f\left( \overline{z_{1}} \right)} = \overline{f\left( \overline{z_{2}} \right)}\]\[\Longrightarrow f\left( \overline{z_{1}} \right) = f\left( \overline{z_{2}} \right)\]

From univalence of \(f(z),\ \overline{z_{1}} = \overline{z_{2}}\). It
follows that \(\overline{f\left( \overline{z} \right)}\ \)is univalent.

\end{proof}

(ii) Rotation: \(e^{- i\theta}f\left( e^{i\theta}z \right) \in S\)

\begin{proof}

Indeed,
\(e^{-i\theta}f\left( e^{i\theta} \cdot 0 \right)=f(0)=0\ \), and

\(\left( e^{- i\theta}f\left( e^{i\theta}z \right) \right)^{'}=e^{- i\theta}e^{i\theta} f'(0)\ = f'(0)=  1\). Univalence follows from \(f(z)\), times a constant. 

\end{proof}

(iii) Dilation: \(r^{- 1}f(r \cdot z\ ) \in S,\ 0 < r < 1\)

\begin{proof}

Indeed, \(r^{- 1}f(r \cdot 0) = f(0) = 0\ \)and
\(\left( r^{- 1}f(r \cdot z\ ) \right)^{'} = rr^{- 1}f'(r \cdot 0)\  = f'(0) = 1\). Univalence follows from \(f(z).\)

\end{proof}

(iv) Disk automorphism

\[\frac{f\left( \frac{z + a}{1 + \overline{a}z} \right) - f(a)}{\left( 1 - |a|^{2} \right) \cdot f'(a)} \in S,\ \ \ |a| < 1\]

\begin{proof}

\[\frac{f\left( \frac{0 + a}{1 + 0} \right) - f(a)}{\left( 1 - |a|^{2} \right) \cdot f'(a)} = \frac{f(a) - f(a)}{\left( 1 - |a|^{2} \right) \cdot f'(a)} = 0\]

Also

\[\left( \frac{f\left( \frac{z + a}{1 + \overline{a}z} \right) - f(a)}{\left( 1 - |a|^{2} \right) \cdot f^{'}(a)} \right)^{'} = \frac{f^{'}\left( \frac{z + a}{1 + \overline{a}z} \right)}{\left( 1 + \overline{a}z \right)^{2} \cdot f^{'}(a)}\]

Substituting\(\ z = 0\) gives\[\frac{f^{'}(a)}{1 \cdot f^{'}(a)} = 1\]

Univalence:

\[\frac{f\left( \frac{z_{1} + a}{1 + \overline{a}z_{1}} \right) - f(a)}{\left( 1 - |a|^{2} \right) \cdot f'(a)} = \frac{f\left( \frac{z_{2} + a}{1 + \overline{a}z_{2}} \right) - f(a)}{\left( 1 - |a|^{2} \right) \cdot f'(a)}\]

\[\Longrightarrow f\left( \frac{z_{1} + a}{1 + \overline{a}z_{1}} \right) = f\left( \frac{z_{2} + a}{1 + \overline{a}z_{2}} \right)\]

From univalence of \(f,\ \)

\[\frac{z_{1} + a}{1 + \overline{a}z_{1}} = \frac{z_{2} + a}{1 + \overline{a}z_{2}}\]

\[\ {\Longrightarrow z}_{1}\left( 1 - |a|^{2} \right) = z_{2}\left( 1 - |a|^{2} \right)\]

\[\Longrightarrow z_{1} = z_{2}\]

\end{proof}

(v) Range transformation.

(vi) Omitted-value transformation.

(vii)Square-root transformation.

The last three results are well-established and proved with details in
\hyperlink{D1}{[D1]}, p.27.

\newpage

For each\(\ f \in S\ \)we can define

\[g(z) = \frac{1}{f\left( z^{- 1} \right)} = z - a_{2} + \left( a_{2}^{2} - a_{3} \right)z^{- 1} + \ldots\]

which is analytic and univalent in the domain
\(\left\{ z:|z| > 1 \right\}.\ \)We are hence motivated to define
another class of univalent functions, class\(\ \Sigma,\ \)which is
analytic and univalent in the domain\(\ \left\{ z:|z| > 1 \right\}:\)

\[g(z) = z + b_{0} + b_{1}z^{- 1} + b_{2}z^{- 2} + \ldots\]

We also define the subclass\(\ \Sigma^{'} \in \Sigma\) for
\(g \in \Sigma^{'}\ \)containing the origin in its image.

The previous transformation, referred as an \emph{inversion},
establishes an one-to-one correspondence
between\(\ S\ \)and\(\ \Sigma^{'}.\ \)For more details, consult \hyperlink{D1}{[D1]}
p.28.

\subsection{1.3 Preliminary results}

We next establish some important elementary results concerning the
classes\(\ S\ \)and\(\ \Sigma.\) We start with a classical and important
theorem in the study of univalent functions.

\uline{Theorem 1.3.1 (the Area Theorem)}

\fbox{\parbox{1.0\textwidth}{\ 
\emph{If}\(\ g \in \Sigma,\ \)\emph{then}

\[\sum_{n = 1}^{\infty}{n\left| b_{n} \right|^{2} \leq 1}\]
}}

\begin{proof}

\emph{(following {[}D1{]})}

Let\(\ E\ \)be the set \textbf{omitted}
by\(\ g.\ \)For\(\ r > 1,\ \)let\(\ C_{r}\ \)be the image
under\(\ g\ \)of the circle\(\ |z|\).

We will proof later that since\(\ g\ \)is univalent, \(\ C_{r}\ \)is a
simple closed curve which encloses a
domain\(\ E_{r},\ \)and\(\ E \subset E_{r}.\ \)By Green's theorem, the
area\(\ A_{r}\ \)of\(\ E_{r}\ \)is given by

\[{\ \ \ A}_{r} = \frac{1}{2i}\int_{C_{r}}^{}(ix\ dy - iy\ dx)\]

\[= \frac{1}{2i}\int_{C_{r}}^{}\overline{w}dw\]

\[= \frac{1}{2i}\int_{|z| = r}^{}{\overline{g(z)}g^{'}(z)}dz\]

\[= \frac{1}{2i}\int_{|z| = r}^{}{\left( \overline{z + b_{0} + b_{1}z^{- 1} + b_{2}z^{- 2} + \ldots} \right)\left( 1 - \sum_{n = 1}^{\infty}{nb_{n}z^{- n - 1}} \right)}dz\]

\[= \frac{1}{2}\int_{0}^{2\pi}{\left\lbrack r^{2} + \sum_{n = 0}^{\infty}{\overline{b_{n}}r^{1 - n}e^{i(n + 1)\theta}} \right\rbrack \cdot \left\lbrack 1 - \sum_{v = 1}^{\infty}{vb_{v}r^{- v - 1}e^{- i(v + 1)\theta}} \right\rbrack}d\theta\]

\[= \frac{1}{2}\left\{ \int_{0}^{2\pi}\left( r^{2} \right)d\theta \right\} + \frac{1}{2}\left\{ \int_{0}^{2\pi}\left( \sum_{n = 0}^{\infty}{\overline{b_{n}}r^{1 - n}e^{i(n + 1)\theta}} \right)d\theta \right\}\]

\[- \frac{1}{2}\left\{ \int_{0}^{2\pi}\left( \sum_{v = 1}^{\infty}{vb_{v}r^{- v + 1}e^{- i(v + 1)\theta}} \right)d\theta \right\} \]

\[- \frac{1}{2}\left\{ \int_{0}^{2\pi}\left( \sum_{n = 0}^{\infty}{\overline{b_{n}}r^{1 - n}e^{i(n + 1)\theta}} \cdot \sum_{v = 1}^{\infty}{vb_{v}r^{- v - 1}e^{- i(v + 1)\theta}} \right)d\theta \right\}\]

\[= \pi r^{2} + \frac{1}{2}\left\{ \sum_{n = 0}^{\infty}{\overline{b_{n}}r^{- n + 1}}\int_{0}^{2\pi}\left( e^{i(n + 1)\theta} \right)d\theta \right\} - \frac{1}{2}\left\{ \sum_{v = 1}^{\infty}{vb_{v}r^{- v + 1}}\int_{0}^{2\pi}\left( e^{- i(v + 1)\theta} \right)d\theta \right\} \]

\[- 
\frac{1}{2}\left\{ \sum_{v = 1}^{\infty}{\sum_{n = 0}^{\infty}\overline{b_{n}} \cdot vb_{v}}\int_{0}^{2\pi}\left( e^{(n - v)i\theta} \right)d\theta \right\}\]

by Cauchy product and taking out the constant. ~\\

Note that if \(k \neq 0\)

\[\int_{0}^{2\pi}\left( e^{- ik\theta} \right) d\theta = \left\lbrack \frac{e^{- ik\theta}}{- ik} \right\rbrack_{0}^{2\pi} = 0 - 0 = 0\]

Hence

\[\int_{0}^{2\pi}\left( e^{- i(m + 1)\theta} \right)d\theta = \int_{0}^{2\pi}\left( e^{i(n + 1)\theta} \right)d\theta = 0\]

Therefore

\[\ \ \ A_{r} = \pi r^{2} + 0 - 0 - \frac{1}{2}\left\{ \sum_{v = 1}^{\infty}{\sum_{n = 0}^{\infty}\overline{b_{n}} \cdot vb_{v}}\int_{0}^{2\pi}\left( e^{(n - v)i\theta} \right)d\theta \right\}\]

\[= \pi r^{2} - \frac{1}{2}\left\{ \sum_{n = 1}^{\infty}{n\left| b_{n} \right|^{2}}\int_{0}^{2\pi}\left( e^{(n - n)i\theta} \right)d\theta \right\}\]

\[= \pi r^{2} - \frac{1}{2}\left\{ \sum_{n = 1}^{\infty}{n\left| b_{n} \right|^{2}}2\pi \right\}\]

\[= \pi\left\{ r^{2} - \sum_{n = 1}^{\infty}{n\left| b_{n} \right|^{2}r^{- 2n}} \right\}\]

The theorem then follows from the fact that,
letting\(\ r \downarrow 1,\ \)

\[A_{r} \rightarrow \pi\left\{ 1 - \sum_{n = 1}^{\infty}{n\left| b_{n} \right|^{2}} \right\} \geq 0\ \iff \sum_{n = 1}^{\infty}{n\left| b_{n} \right|^{2}} \leq 1\]

since the area of\(\ A_{r} \geq 0.\ \)

\end{proof}

~\\

In the above proof, we made an inportant assumption:

\hypertarget{Int}{\uline{Proposition 1.3.2}}

\fbox{\parbox{1.0\textwidth}{\ 
Consider\(\ g\ \)on the domain
\(\mathbb{C}^{\infty}\backslash\left\{ D(0,1) \right\}.\ \)Then every
point in the interior
of\(\ \mathbb{C}^{\infty}\backslash\left\{ D(0,r) \right\}\ \)is mapped
to the interior
of\(\ \ \mathbb{C}^{\infty}\backslash\left\{ E_{r} \right\}\), and moreover, we have that \(\ E \subset E_{r},\ \)where\(\ E\ \)is as defined above.
}}

Although this seems right intuitively, it is not as easy as it seems to
be proved, and leads doubt and questions to reader's view on the
validity of the above proof. We hence present a proof of this statement.
The proof is in two phases:

\uline{Lemma 1.3.3}

\fbox{\parbox{1.0\textwidth}{\ 
\emph{Suppose}\(\ f(z)\ \)\emph{is an univalent holomorphic function,
mapping a simply connected domain}\(\mathbb{\ D}\) \emph{to some
domain}\(\mathbb{\ W}\)\emph{. Then for an arbitrary Jordan
curve}\(\ C \subset \mathbb{D\ }\)\emph{mapped under}\(\ f\ \)\emph{onto
a Jordan curve}\(\ \Gamma \subset \mathbb{W,\ }int(C)\ \)\emph{is mapped
1-to-1 to}\(\ int(\Gamma)\).
}}
\newpage
\begin{proof}

Let\(\ w_{0}\ \)be an arbitrary point in\(\mathbb{\ W\ }\)and consider
the set of points\(\ z \in int(C)\). We consider the number of times
that\(\ f(z)\ \)could take\(\ z\ \)to\(\ w_{0}\ \)by counting the number
of zeroes obtained by \(g(z) = f(z) - w_{0}\). 

That is,

\[N - P = \frac{1}{2\pi i}\int_{C}^{}{\frac{g^{'}(z)}{g(z)}dz}\]

where\(\ N\ \)is the number of zeroes
of\(\ g\ \)inside\(\ C,\ \)and\(\ P\ \)is the number of poles
of\(\ g\ \)inside\(\ C,\ \) positively oriented in the integral.

Notice that\(\ P = 0\ \)as\(\ f\ \)is holomorphic in\(\ C\), hence

\[N = \frac{1}{2\pi i}\int_{C}^{}{\frac{g^{'}(z)}{g(z)}dz} = \frac{1}{2\pi i}\int_{C}^{}{\frac{f^{'}(z)}{f(z) - w_{0}}dz}\]

Parameterizing\(\ C\ \)and\(\ \Gamma\ \)respectively, we have that

\[\left\{ \begin{matrix}
C = \gamma(t),\ \ t \in \lbrack 0,a\rbrack \\
\Gamma = f\left( \gamma(t) \right),\ \ \ t \in \lbrack 0,a\rbrack \\
\end{matrix} \right.\ \]

Notice that

\[\int_{\Gamma}^{}\frac{1}{z - w_{0}}dz = \int_{0}^{a}\frac{f^{'}\left( \gamma(t) \right)\gamma^{'}(t)}{f\left( \gamma(t) \right) - w_{0}}dt\]

by the chain rule. It follows that
\begin{align*}
N &=\frac{1}{2\pi i}\int_{C}^{}{\frac{f^{'}(z)}{f(z) - w_{0}}dz}\\
&=\frac{1}{2\pi i}\int_{0}^{a}\frac{f^{'}\left( \gamma(t) \right)}{f\left( \gamma(t) \right){- w}_{0}}\gamma^{'}(t)dt\\
&=\frac{1}{2\pi i}\int_{0}^{a}\frac{f^{'}\left( \gamma(t) \right)\gamma^{'}(t)}{f\left( \gamma(t) \right) - w_{0}}dt\\
&=\frac{1}{2\pi i}\int_{\Gamma}^{}\frac{1}{z - w_{0}}dz
\end{align*}

So that the number of times that\(\ f(z)\ \)could
take\(\ z\ \)to\(\ w_{0}\ \)is now dependent on an integral on the
positively oriented contour\(\ \Gamma\ \)instead of\(\ C\).

\uline{Case I.}

If\(\ w_{0}\ \)lies outside of\(\ int(\Gamma)\), then

\[N = \int_{\Gamma}^{}{\frac{1}{w - w_{0}}dw} = 0\]

by Cauchy's theorem since\(\ \Gamma\ \)is closed.

~\\
Then\(\ N = 0\ \)implies that\(\ f\ \)could never take some \(z \in int(C)\ \)to\(\ w_{0}\ \)that is outside
of\(\ int(\Gamma)\). 

\uline{Case II.}

If\(\ \omega\ \)lies inside\(\ \Gamma\), then

\[\int_{\Gamma}^{}{\frac{1}{w - w_{0}}dw} = 2\pi i \Longrightarrow N = 1\]

Only this case is valid,
hence\(\ f(z)\ \)maps\(\ int(C)\ \)to\(\ int(\Gamma)\) injectively by
the univalence of\(\ f.\)

Since\(\ w_{0}\ \)is arbitrary, this implies that every point inside the
interior of\(\ \Gamma\ \)has a pre-image from the interior
of\(\ C.\ \)This conclusion is useful in the next given proof.

\end{proof}

\uline{Corollary 1.3.4}

\fbox{\parbox{1.0\textwidth}{\ 
\emph{Suppose}\(\ f(z)\ \)\emph{is an univalent holomorphic function
defined on the
domain}\(\mathbb{\ D.\ }\)\emph{If}\(\mathbb{\ D\ }\)\emph{is simply
connected then so is its
image}\(\ \Delta = f\left( \mathbb{D} \right)\ \)\emph{.}
}}

\begin{proof}

For an arbitrary Jordan curve\(\ \Gamma \subset \Delta,\ \)take an
arbitrary inscribed polygon\(\ P \subset \Gamma.\ \)Then\(\ P\ \)is the
image of some simply closed curve\(\ C\ \)in\(\mathbb{\ D.\ }\)Then by
Lemma 1.3.3, the interior of\(\ C\ \)is mapped injectively onto the
interior of\(\ P\). Since any point in the interior of\(\ P\ \)is a
point in the interior of\(\ \Gamma\ \)hence also the interior of
\(\Delta,\) \(\Delta\ \)is simply connected by definition
since\(\ P\ \)is arbitrary.

\end{proof} 

Alternatively, one could use homotopy to proof the same result. Readers
who are not familiar with homotopy are encouraged to skip to the next
proposition.

\uline{Lemma 1.3.3(Via homotopy)}

\fbox{\parbox{1.0\textwidth}{\ 
\emph{Suppose}\(\ f(z)\ \)\emph{is an univalent holomorphic function,
mapping a simply connected domain}\(\mathbb{\ D}\) \emph{to some
domain}\(\mathbb{\ W}\)\emph{. Then for an arbitrary Jordan
curve}\(\ C \subset \mathbb{D\ }\)\emph{mapped under}\(\ f\ \)\emph{onto
a Jordan curve}\(\ \Gamma \subset \mathbb{W,\ }int(C)\ \)\emph{is mapped
1-to-1 to}\(\ int(\Gamma)\).
}}

\begin{proof}

Let\(\ \gamma_{1}\ \)be an arbitrary closed curve in\(\ G_{1}\).
Since\({\ G}_{1}\ \)is simply connected, \(\gamma_{1}\ \)is homotopic to
zero, hence there is a continuous map

\[\phi:\lbrack 0,1\rbrack \times \lbrack 0,1\rbrack \rightarrow {\ G}_{1}\]

where

\[\left\{ \begin{matrix}
\phi(t,0) = \gamma_{1}(t)\ \ \ \ \ \ \ \ \ \ \ \ \ \ \ \ \ \ \ \ \ \ \ \ \ \ \ t \in \lbrack 0,1\rbrack \\
\phi(t,1) = z_{0}\ \ \ \ \ \ \ \ \ \ \ \ \ \ \ \ \ \ \ \ \ \ \ \ \ \ \ \ \ \ \ \ \ t \in \lbrack 0,1\rbrack \\
\phi(0,s) = \phi(1,s)\ \ \ \ \ \ \ \ \ \ \ \ \ \ \ \ \ \ \ \ \ \ \ s \in \lbrack 0,1\rbrack \\
\end{matrix} \right.\ \]

for some \(z_{0} \in {\ G}_{1}\)

Suppose the two domains are conformally equivalent under some continuous
analytic function\(\ f\). Then
\(f(\phi):\lbrack 0,1\rbrack \times \lbrack 0,1\rbrack \rightarrow {\ G}_{2}\)
is a composition of continuous functions, hence is also continuous.
Since the there exists an inverse function for all conformal mappings,
\({\ \gamma}_{1} = f^{- 1}\left( \gamma_{2}(t) \right)\) for some curve
in\({\ \ G}_{2},\ \)where\({\ f}^{- 1}\left( f(z) \right)\) is the
identity mapping. Note that

\[f\left( \phi(t,0) \right) = f\left( \gamma_{1}(t) \right) = f\left( f^{- 1}\left( \gamma_{2}(t) \right) \right) = \gamma_{2}(t)\]

\begin{quote}
\(f\left( \phi(t,1) \right) = f\left( z_{0} \right),\ \)which is a point
in\(\ {\ G}_{2}\)

\(f\left( \phi(0,s) \right) = f\left( \phi(1,s)\  \right),\ \)from the
univalence of\(\ f\)
\end{quote}

Hence it can be concluded that the continuous
map\(\ f(\phi)\ \)constructs a path homotopy from
\(\gamma_{2}\ \)to\(\ f\left( z_{0} \right)\).
Since\({\ \gamma}_{2}\ \)is an arbitrary curve in\({\ \ G}_{2}\ \)by
choosing a corresponding\({\ \gamma}_{1},\ \)this implies
that\({\ \ G}_{2}\ \)is indeed simply connected.

\end{proof}

\uline{(\hyperlink{Int}{Proposition 1.3.2} revisited)}

\fbox{\parbox{1.0\textwidth}{\ 
\(E \subset E_{r}\ \)for the class of univalent
functions\(\ g \in \Sigma\).
}}

\begin{proof}

Consider the domain
of\(\ g \in \Sigma,\ \ C^{\infty}\backslash\left\{ \mathbb{D} \right\},\ \)stereographically
projected onto the Riemann Sphere. Then since every closed curve is
homologically equivalent to a point,
\(C^{\infty}\backslash\left\{ \mathbb{D} \right\}\ \)is a simply
connected domain. By Corollary 1.3.4, the image of this domain
under\(\ g\ \)must also be simply connected.

Consider\(\ g\ \)on the domain
\(C^{\infty}\backslash\left\{ D(0,r) \right\}.\ \)Then by Lemma 1.3.3,
every point in the interior
of\(\ C^{\infty}\backslash\left\{ D(0,r) \right\}\ \)is mapped to the
interior of\(\ \ C^{\infty}\backslash\left\{ E_{r} \right\}\). Suppose
for the sake of contradiction that\(\ E\ \)is outside\(\ E_{r},\ \)then
the image is no longer a simply connected domain, contradicting to
Corollary 1.3.4. Hence\(\ E \subset E_{r}.\ \)This finishes our
discussion on Proposition 1.3.2, and we are now convinced that indeed,
we can let\(\ r \rightarrow 1\ \)in Theorem 1.3.1.

 \end{proof}

The following corollaries follows from the Area theorem; they can be
found on {[}D1{]}, p.30 for their statements without proves. Here we
give proves to them:

\uline{Corollary 1.3.5}

\fbox{\parbox{1.0\textwidth}{\ 
If\(\ g \in \Sigma,\ \)then

(i)

\[\left| b_{n} \right| \leq n^{\frac{- 1}{2}}\ \ \ \forall n\mathbb{\in N}\]

(ii)

\[\left| b_{1} \right| \leq 1\]

\[\left| b_{1} \right| = 1\iff g(z) = z + b_{0} + \frac{b_{1}}{z}\ \ \]

(iii)

This is a conformal mapping of \(\Delta\ \)onto the complement of a line
segment with length 4.
}}

\begin{proof}

(i) It follows from the Area Theorem (Theorem 1.3.1) that

\[\sum_{n = 1}^{\infty}{n\left| b_{n} \right|^{2}} \leq 1 \Longrightarrow n\left| b_{n} \right|^{2} \leq 1\ \]

fir any arbitrary\(\ n\ \)

\[\Longrightarrow \left| b_{n} \right|^{2} \leq n^{- 1}\]

\[\Longrightarrow \left| b_{n} \right| \leq n^{- \frac{1}{2}}\ \ \]

(ii)

From (i), \(\left| b_{1} \right| \leq 1\).
Suppose\(\ \left| b_{1} \right| = 1\), then since

\[\sum_{n = 1}^{\infty}{n\left| b_{n} \right|^{2}} \leq 1\]

from the Area theorem. Note that \(\left| b_{n} \right| = 0,\ \ \forall n \geq 2\), or else

\[\sum_{n = 1}^{\infty}{n\left| b_{n} \right|^{2}} > 1\]

as \(\left| b_{1} \right| = 1\ \)already. It follows that
\(g(z) = z + b_{0} + \frac{b_{1}}{z}\ \)as later coefficients are zero.

(iii)

To see that the codomain of\(\ g\ \)is a line segment with length 4,
first set \(|z| = 1\). Then

\[g(z) = z + b_{0} + b_{1} \cdot \overline{z}\]

Rearrange\(\ g\ \)to get

\[\ \ \ \ g(z) = z + b_{0} + b_{1} \cdot \overline{z} + b_{1} \cdot z - b_{1} \cdot z\]

\[= b_{1} \cdot \left( z + \overline{z} \right) + b_{0} + \left( 1 - b_{1} \right) \cdot z\]

\[= b_{1} \cdot 2\mathfrak{R}(z) + b_{0} + \left( 1 - b_{1} \right) \cdot z\]

Notice that\(\ b_{1} \cdot 2\mathfrak{R}(z)\ \)transforms
all\(\ z\ \)to a line segment with endpoint 2\(\ b_{1}\)and
-2\(\ b_{1},\ \) and has a length
of\(\ 4 \cdot \left| b_{1} \right| = 4 \cdot 1 = 4\), as the value
of \(\mathfrak{R}(z)\ \)lies between 1 and \(- 1\); \(b_{0}\ \)is a
complex number that transforms every single point on this line segment
to a certain direction without changing the length, and so does
\(\left( 1 - b_{1} \right) \cdot z;\ \)it then follows that setting
\(|z| > 1\) then the image of\(\ g\ \)avoids this line segment, hence is
a conformal mapping of \(\Delta\ \)onto the complement of a line segment
with length 4.

\end{proof} 

We also present two other frequently used, well-known results. Since
these results have detailed proves in {[}D1{]} p.30-31, we only state
them for later uses.

\uline{Theorem 1.3.6 (Bieberbach's Theorem)}

\fbox{\parbox{1.0\textwidth}{\ 
\emph{If}\(\ f \in S,\ \)\emph{then}\(\ \left| a_{n} \right| \leq 2,\ \)\emph{with
equality iff}\(\ f\ \)\emph{is a rotation of the Koebe Function.}
}}

\hypertarget{Koebe}{\uline{Theorem 1.3.7 (Koebe One-Quarter Theorem)}}

\fbox{\parbox{1.0\textwidth}{\ 
\emph{If}\(\ f \in S,\ \)\emph{then}\(\ f\left( \mathbb{D} \right)\)
\emph{always contains the set}\(\ \left\{ w:|w| < 1/4 \right\}.\)}}

The above two theorems are interesting in nature: Theorem 1.3.6 makes a
start to the famous Bieberbach's conjecture (conjecture 2.1.1), whilst
Theorem 1.3.7 asserts that however we play around with a function in
class\(\ S,\ \)we always end up having a disc with radius\(\ 1/4\)
within its image. The later theorem is particularly useful in later
arguments.

\subsection{1.4 Further Theorems}

We now present some further theorems on univalent function theory. Since
these results have been proven in {[}D1{]} p.32-40, we only make
supplements to them when necessary.

\uline{Theorem 1.4.1}

\fbox{\parbox{1.0\textwidth}{\  
\emph{If}\(\ f \in S,\ \)\emph{then}

\[\left| \frac{zf^{''}(z)}{f'(z)} - \frac{2r^{2}}{1 - r^{2}} \right| \leq \frac{4r}{1 - r^{2}}\]

\(|z| = r < 1\).
}}

The proof in {[}D1, p.32{]} omits the explanation in the computation
of\(\ A_{2\ }\); here we complete these parts:

\begin{proof}

Given\(\ f \in S,\ \)fix\(\ \xi \in \mathbb{D\ }\)and by a disk
automorphism (recall from chapter 1.2 transformation (iv) that this
preserves\(\ f \in S\)) we define

\[F(z) = \frac{f\left( \frac{z + \xi}{1 + \overline{\xi}z} \right) - f(\xi)}{\left( 1 - |\xi|^{2} \right)f^{'}(\xi)}\]

Note that since\(\ F\ \)has the Taylor expansion

\[F(z) = z + A_{2}z^{2} + A_{3}z^{3} + \ldots\]

\[{\Longrightarrow F}^{''}(z) = 2A_{2} + 6A_{3}z + \ldots\]

\[\Longrightarrow F^{''}(0) = 2A_{2}\]

Hence if follows from

\[F^{''}(0) = \left( 1 - |\xi|^{2} \right)\frac{f^{''}(\xi)}{f'(\xi)} - 2\overline{\xi}\]

that

\[A_{2} = \frac{1}{2}F^{''}(0) = \frac{1}{2}\left\lbrack \left( 1 - |\xi|^{2} \right)\frac{f^{''}(\xi)}{f'(\xi)} - 2\overline{\xi} \right\rbrack\]

By Bieberbach's Theorem

\[\left| A_{2} \right| = \left| \left( 1 - |\xi|^{2} \right)\frac{f^{''}(\xi)}{2f'(\xi)} - \overline{\xi} \right| \leq 2\]

Multiplying both sides by

\[\frac{2|\xi|}{1 - |\xi|^{2}}\]

gives

\[\left| \frac{\xi f^{''}(\xi)}{f'(\xi)} - \frac{2|\xi|^{2}}{1 - |\xi|^{2}} \right| \leq \frac{4|\xi|}{1 - |\xi|^{2}}\]

which is equivalent to the desired form

\[\left| \frac{zf^{''}(z)}{f'(z)} - \frac{2r^{2}}{1 - r^{2}} \right| \leq \frac{4r}{1 - r^{2}}\]

 \end{proof}

\hypertarget{Dist}{\uline{Theorem 1.4.2 (the Distortion Theorem)}}

\fbox{\parbox{1.0\textwidth}{\ 
\emph{If}\(\ f \in S,\ \)\emph{then}

\[\frac{1 - r}{(1 + r)^{3}} \leq \left| f^{'}(z) \right| \leq \frac{1 + r}{(1 - r)^{3}}\]

\(|z| = r < 1\).

}}

\begin{proof}
    
 See \hyperlink{D1}{[D1]} p.32. In \hyperlink{D1}{[D1]} it has been shown that

\[\frac{2r - 4}{1 - r^{2}} \leq \frac{\partial}{\partial r}\log{\left| f^{'}\left( re^{i\theta} \right) \right| \leq}\frac{2r + 4}{1 - r^{2}}\]

We integrate, which preserves the inequality, from\(\ 0\ \)to\(\ R,\ \)

\[\int_{0}^{R}\frac{2r - 4}{1 - r^{2}}dr \leq \int_{0}^{R}\left( \frac{\partial}{\partial r}\log\left| f^{'}\left( re^{i\theta} \right) \right| \right)dr \leq \int_{0}^{R}\frac{2r + 4}{1 - r^{2}}dr\]

as\(\ r < 1\ \)and all functions are analytic

\[\log\frac{1 - R}{(1 + R)^{3}} \leq \int_{0}^{R}\left( \frac{\partial}{\partial r}\log\left| f^{'}\left( re^{i\theta} \right) \right| \right)dr \leq \log\frac{1 + R}{(1 + R)^{3}}\]

For

\[\int_{0}^{R}{\left( \frac{\partial}{\partial r}\log\left| f^{'}\left( re^{i\theta} \right) \right| \right)dr}\]

By the fundamental theorem of calculus

\[\int_{0}^{R}{\left( \frac{\partial}{\partial r}\log\left| f^{'}\left( re^{i\theta} \right) \right| \right)dr =}\log\left| f^{'}\left( Re^{i\theta} \right) \right|\]

and hence the theorem follows.

 \end{proof}

\hypertarget{Growth}{\uline{Theorem 1.4.3 (the Growth Theorem)}}

\fbox{\parbox{1.0\textwidth}{\ 
\emph{If}\(\ f \in S,\ \)\emph{then}

\[\frac{r}{(1 + r)^{2}} \leq \left| f(z) \right| \leq \frac{r}{(1 - r)^{2}}\]

\(|z| = r < 1\).
}}

\begin{proof}
    
 See \hyperlink{D1}{[D1]} p.34.\end{proof}

Since\(\ \left| f(z) \right|\ \)is bounded above
by\(\ r/(1 - r)^{2},\ \)we have by Montel's theorem (See \hyperlink{Montel}{Thm 1.1.5})

\uline{Theorem 1.4.4}

\fbox{\parbox{1.0\textwidth}{\ 
\emph{The class}\(\ S\ \)\emph{of univalent functions is a compact
normal family.}

}}

For further explanations consult \hyperlink{D1}{[D1]]} p.9. We finish this chapter by
presenting a theorem that is a direct consequence from \hyperlink{Growth}{Theorem 1.4.3}.

\uline{Theorem 1.4.5}

\fbox{\parbox{1.0\textwidth}{\ 

\emph{If}\(\ f \in S,\ \)\emph{then}

\[\frac{1 - r}{1 + r} \leq \left| \frac{zf^{'}(z)}{f(z)} \right| \leq \frac{1 + r}{1 - r}\]

\(|z| = r < 1\).

}}

\emph{Proof.} See \hyperlink{D1}{[D1]]} p.35

\newpage

\hypertarget{the-bieberbach-conjecture-1}{%
\section{2 The Bieberbach
conjecture}\label{the-bieberbach-conjecture-1}}

\subsection{2.1 The Bieberbach conjecture}

As many might have guessed, the Bieberbach conjecture concerns putting
an upper bound on the coefficients of the series
representing\(\ f \in S\). The proposed upper bound itself is remarkably
straight forward:

\uline{Theorem 2.1.1 (the Bieberbach Conjecture/de Brange's theorem)}

\fbox{\parbox{1.0\textwidth}{\ 
\emph{The coefficients of each function}\(\ f \in S\ \)\emph{satisfies}\(\ \left| a_{n} \right| \leq n\ \ \)\emph{for}\(\ \ n = 2,3,\ldots\)

\emph{Equality occurs iff the function}\(\ f\ \)\emph{is a rotation of
the Koebe Function or its rotations.}

}}

The Bieberbach Conjecture was proposed by Ludwig Bieberbach in 1916 and
finally proven by Louis de Branges in 1985. Over the course of 70 years,
many techniques have been developed in attempt to prove this seemingly
easy proposition, some leading to very modern developments in the field
of Complex Analysis, still actively being studied today.

To familiarize ourselves with how mathematicians attempted this famous
conjecture at the very start, we introduce one important attempt to
proof the Bieberbach conjecture, known as Littlewood's theorem, proved
in 1925.

\hypertarget{Thm2.1.2}{\uline{Theorem 2.1.2 (Littlewood's theorem)}}

\fbox{\parbox{1.0\textwidth}{\ 
\emph{The coefficients of each
family}\(\ f \in S\ \)\emph{satisfies}\(\ \left| a_{n} \right| \leq e \times n\ \)\emph{for}\(\ n = 2,3,\ldots\)

}}

\begin{proof}

We define for\(\ f \in S\)

\[M_{p}(r,f) = \left\{ \frac{1}{2\pi}\int_{0}^{2\pi}\left| f\left( re^{i\theta} \right) \right|^{p}d\theta \right\}^{\frac{1}{p}},\ 0 < p < \infty\]

for\(\ 0 < p < \infty.\ \)We first proof that

\[M_{p}(r,f) \leq \frac{r}{1 - r}\]

Given\(\ f \in S,\ \)let

\[h(z) = \sqrt{f\left( z^{2} \right)} = \sum_{n = 1}^{\infty}{c_{n}z^{n}}\]

which in view of the Elementary Transformations (Recall from \hyperlink{Ele}{Chapter 1.2})
still belongs to the class\(\ S.\ \)Hence the growth theorem (See \hyperlink{Growth}{Theorem 1.4.3}) applies, so

\[\left| f(z) \right| \leq \frac{r}{(1 - r)^{2}} \Longrightarrow f\left( z^{2} \right) \leq \frac{r^{2}}{\left( 1 - r^{2} \right)^{2}} \Longrightarrow \sqrt{f\left( z^{2} \right)} = h(z) \leq \frac{r}{1 - r^{2}}\]

In view of Lemma 1.3.3, \(\ h\ \)maps the disc\(\ |z| < r\) conformally
onto a domain\(\ D_{r}\ \)that lies within the
disk\(\ |w| < \frac{r}{1 - r^{2}}.\ \) The area\(\ A_{r}\)
of\(\ D_{r}\ \)is therefore no greater than this disk:

\[A_{r} \leq \pi\left( \frac{r}{1 - r^{2}} \right)^{2}\]

where we used the standard formula\(\ A = \pi r^{2}\) for the area of a
circle. We may also explicitly calculate\(\ A_{r}\ \)giving

\[A_{r} = \int_{0}^{2\pi}{\int_{0}^{r}{\left| h^{'}\left( \rho e^{i\theta} \right) \right|^{2}\rho d\rho d\theta}} = \int_{0}^{2\pi}{\int_{0}^{r}{\sum_{n = 1}^{\infty}{n^{2}\left| c_{n} \right|^{2}}\rho^{2n - 1}d\rho d\theta}}\]

By differentiating and squaring \(h(z)\ \)using series definition. This simplifies to

\[A_{r} = \int_{0}^{2\pi}{\frac{1}{2n}\sum_{n = 1}^{\infty}{n^{2}\left| c_{n} \right|^{2}}r^{2n}d\theta} = \pi\sum_{n = 1}^{\infty}{n\left| c_{n} \right|^{2}}r^{2n}\]

Consequently

\[\sum_{n = 1}^{\infty}{n\left| c_{n} \right|^{2}}r^{2n} \leq \frac{r^{2}}{\left( 1 - r^{2} \right)^{2}} \Longrightarrow \sum_{n = 1}^{\infty}{n\left| c_{n} \right|^{2}}r^{2n - 1} \leq \frac{r}{\left( 1 - r^{2} \right)^{2}} \Longrightarrow \sum_{n = 1}^{\infty}\left| c_{n} \right|^{2}r^{2n} \leq \frac{r^{2}}{1 - r^{2}}\]

where we first divided by\(\ r\ \)then integrated term-by-term
from\(\ 0\ \)to\(\ r\ \)on both sides. This in turn shows that

\[\frac{1}{2\pi}\int_{0}^{2\pi}\left| h\left( re^{i\theta} \right) \right|^{2}d\theta \leq \frac{r^{2}}{1 - r^{2}}\]

From the above argument we have that indeed

\[M_{p}(r,f) \leq \frac{r}{1 - r}\]

From this step, we note that

\[\ \ \ a_{n} = \frac{1}{2\pi i}\oint_{|z| = r}^{}{\frac{f(z)\ }{z^{n + 1}}dz\ }\]

\[\ \ \  \Longrightarrow \left| a_{n} \right| = \left| \frac{1}{2\pi i}\oint_{|z| = r}^{}{\frac{f(z)\ }{z^{n + 1}}dz\ } \right|\]

\[= \frac{1}{2\pi}\left| \int_{0}^{2\pi}{\frac{f\left( re^{i\theta} \right)}{r^{n}e^{in\theta}}d\theta} \right|\]

\[\leq \frac{1}{2\pi}\int_{0}^{2\pi}{\left| \frac{f\left( re^{i\theta} \right)}{r^{n}e^{in\theta}} \right|d\theta}\]

\[= \frac{1}{2\pi r^{n}}\int_{0}^{2\pi}{\left| f\left( re^{i\theta} \right) \right|d\theta}\]

\[= \frac{1}{r^{n}} \cdot \left\{ \frac{1}{2\pi}\int_{0}^{2\pi}\left| f\left( re^{i\theta} \right) \right|d\theta \right\}\]

\[= {\frac{1}{r^{n}}M}_{1}(r,f) \leq \frac{1}{r^{n}} \cdot \frac{r}{1 - r} = \frac{1}{(1 - r)r^{n - 1}}\]

Which attains a maximum at\(\ r = 1 - 1/n\). Choosing this value
for\(\ r\ \)we have that

\[\left| a_{n} \right| \leq n\left( 1 + \frac{1}{n - 1} \right)^{n - 1} < e \times n\]

completing the proof.

\end{proof}

This is a good approximation of the coefficients\(\ a_{n},\ \)however it
is far from sharp. To estimate\(\ a_{n}\ \)better we need stronger
tools, as we will see in Chapter 3; we also need more subtle arguments,
for instance equivalent statements of the Bieberbach conjecture that are
easier to work with. We spend the next two chapters introducing two of
such equivalent conjectures that played a crucial role in proving the
result that we are looking for.

\subsection{2.2 Robertson's conjecture}

In 1936, M.S.Robertson proposed a conjecture that implies Bieberbach's
conjecture, known as Robertson's conjecture. Later in 1964, I.M.Milin
proposed another conjecture that implies Robertson's conjecture (and in
turn, of course, Bieberbach's conjecture), which was proved by de
Branges in 1985 (and later by Weinstein in 1991).

We spend the next two Chapters introducing the above two conjectures,
and prove that they indeed imply Bieberbach's conjecture.

It has been verified (see \hyperlink{Ele}{Chapter 1.2}) that

\[g(z) = \sqrt{f\left( z^{2} \right)}\]

belongs to the family \(\ S.\ \)

Then

\[g(z) = z + c_{3}z^{3} + c_{5}z^{5} + \ldots\]

\hypertarget{Rob}{\uline{Theorem 2.2.1 (Robertson's conjecture)}}

\fbox{\parbox{1.0\textwidth}{\ 

\[\sum_{k = 1}^{n}{\left| c_{2k - 1} \right|^{2} \leq n}\]

\emph{This conjecture implies Bieberbach's conjecture.}

}}

\begin{proof}

To see this, compare the coefficients of

\[f\left( z^{2} \right) = z^{2} + a_{2}z^{4} + a_{3}z^{6} + \ldots\]

and

\[g(z)^{2} = z^{2} + \left( c_{1}c_{3} + c_{3}c_{1} \right)z^{4} + \left( c_{1}c_{5} + c_{3}c_{3} + c_{5}c_{1} \right)z^{6} + \ldots\]

Notice that

\[a_{n} = c_{1}c_{2n - 1} + c_{3}c_{2n - 3} + \ldots + c_{2n - 1}c_{1}\]

and hence by the triangle inequality

\[\left| a_{n} \right| \leq \left| c_{1} \right|\left| c_{2n - 1} \right| + \left| c_{3} \right|\left| c_{2n - 3} \right| + \ldots + \left| c_{2n - 1} \right|\left| c_{1} \right|\]

It then follows from the Cauchy-Schwartz inequality that

\[\left| a_{n} \right| \leq \left| c_{1} \right|^{2} + \left| c_{3} \right|^{2} + \ldots + \left| c_{2n - 1} \right|^{2} \leq n\]

which is Bieberbach's Conjecture.

\end{proof}

\subsection{2.3 Milin's conjecture}

Let

\[\varphi(z) = \sum_{k = 1}^{\infty}{\alpha_{k}z^{k}}\]

be an arbitrary power series with a positive radius of convergence. We
denote

\[\psi(z) = e^{\varphi(z)} = \sum_{k = 0}^{\infty}{\beta_{k}z^{k}}\]

\hypertarget{LM}{\uline{Theorem 2.3.1 (The Lebedev-Milin Inequality)}}

\fbox{\parbox{1.0\textwidth}{\ 
\emph{For}\(\ n = 1,2,\ldots,\)

\[\sum_{k = 0}^{n}{\left| \beta_{k} \right|^{2} \leq (n + 1)\exp\left\{ \frac{1}{n + 1}\sum_{m = 1}^{n}{\sum_{k = 1}^{m}\left( k\left| \alpha_{k} \right|^{2} - \frac{1}{k} \right)} \right\}}\]
~\\
\emph{Equality occurs for a
given n iff} \[\alpha_{k} = \frac{\gamma^{k}}{k},\ k = 1,2,\ldots,n\ \]

\emph{Here}\(\ \gamma \in \mathbb{C}\)\emph{ is a constant satisfying}\(\ |\gamma| = 1\)\emph{.}

}}

\begin{proof}

Note that\(\ \)

\[\psi(0) = e^{\varphi(0)} = e^{0} = 1 \Longrightarrow \sum_{k = 0}^{\infty}{\beta_{k} \cdot 0^{k}} = 1 \Longrightarrow \beta_{0} = 1\]

To find other\(\ \beta_{n\ }\)s,
differentiating\(\ \psi(z)\ \)gives\(\ \psi^{'}(z) = \psi(z)\varphi^{'}(z).\ \)Differentiating
the power series term by term yields

\[\psi^{'}(z) = \sum_{k = 1}^{\infty}{k\beta_{k}z^{k - 1}} = \psi(z)\varphi^{'}(z) = \left( \sum_{k = 0}^{\infty}{\beta_{k}z^{k}} \right)\left( \sum_{k = 1}^{\infty}{k\alpha_{k}z^{k - 1}} \right)\]

\[\Longrightarrow \sum_{k = 1}^{\infty}{k\beta_{k}z^{k - 1}} = \left( \sum_{k = 0}^{\infty}{\beta_{k}z^{k}} \right)\left( \sum_{k = 1}^{\infty}{k\alpha_{k}z^{k - 1}} \right)\]

comparing coefficients of\(\ z^{n - 1}\ \)gives

\[n\beta_{n}z^{n - 1} = \left( \sum_{k = 0}^{n - 1}(n - k)\alpha_{n - k}\beta_{k} \right)z^{n - 1}\]

Rearranging gives

\[n\beta_{n} = \sum_{k = 0}^{n - 1}(n - k)\alpha_{n - k}\beta_{k}\]

\[\ \ \ \ \ \ \ \  \Longrightarrow \left| n\beta_{n} \right| \leq \sum_{k = 0}^{n - 1}\left| (n - k)\alpha_{n - k} \right|\left| \alpha_{n - k}\beta_{k} \right|\]

\[\Longrightarrow n^{2}\left| \beta_{n} \right|^{2} \leq \sum_{k = 1}^{n - 1}{k^{2}\left| \alpha_{n} \right|^{2}}\sum_{k = 0}^{n - 1}\left| \beta_{n} \right|^{2}\]

by the Cauchy-Schwartz inequality.

Denoting

\[P_{n} = \sum_{k = 1}^{n}{k^{2}\left| \alpha_{k} \right|^{2}}\]

and

\[Q_{n} = \sum_{k = 0}^{n}\left| \beta_{n} \right|^{2}\]

Then

\[\ \ \ Q_{n} = Q_{n - 1} + \left| \beta_{n} \right|^{2}\]

\[\leq Q_{n - 1} + \frac{1}{n^{2}}P_{n}Q_{n - 1}\]

\[= \left\{ 1 + \frac{1}{n^{2}}P_{n} \right\} Q_{n - 1}\]

\[= \frac{n + 1}{n}\left\{ 1 + \frac{P_{n} - n}{n(n + 1)} \right\} Q_{n - 1}\]

\[\leq \frac{n + 1}{n}\exp\left\{ \frac{P_{n} - n}{n(n + 1)} \right\} Q_{n - 1}\]

Note that

\[Q_{0} = \left| \beta_{0} \right|^{2} = 1\]

It follows that

\[Q_{1} \leq 1\exp\left\{ \frac{P_{n} - 1}{2} \right\} Q_{0}\]

\[= 1\exp\left\{ \frac{P_{n} - 1}{2} \right\}\]

\[Q_{2} \leq \frac{3}{2}\exp\left\{ \frac{P_{n} - 2}{6} \right\} Q_{1}\]

\[\ \ \ \ \ \ \ \ \ \ \ \ \ \ \ \ \ \ \ \ \ \ \ \ \ \ \ \ \ \ \ \ \ \ \ \ \ \ \ \  = \left( \frac{3}{2}\exp\left\{ \frac{P_{n} - 2}{6} \right\} \right) \cdot \left( 1\exp\left\{ \frac{P_{n} - 1}{2} \right\} \right)\]

\[\vdots\]

\[Q_{n} \leq \left( 1 \cdot \frac{3}{2} \cdot \frac{3}{4} \cdot \ldots \cdot \frac{n + 1}{n} \right)\exp\left\{ \frac{P_{n} - 1}{2} + \frac{P_{n} - 2}{6} + \ldots\frac{P_{n} - n}{n(n + 1)} \right\}\]

\[= (n + 1)\exp\left\{ \sum_{k = 1}^{n}\frac{P_{k} - k}{k(k + 1)} \right\}\]

\[= (n + 1)\exp\left\{ \sum_{k = 1}^{n}\frac{P_{k}}{k(k + 1)} - \sum_{k = 1}^{n}\frac{1}{k + 1} \right\}\]

\[= (n + 1)\exp\left\{ \sum_{k = 1}^{n}{\frac{P_{k}}{k(k + 1)} + 1 - \sum_{k = 1}^{n + 1}\frac{1}{k}} \right\}\]

Now

\[s_{n} = \sum_{k = 1}^{n}\frac{1}{k(k + 1)} = \sum_{k = 1}^{n}\left( \frac{1}{k} - \frac{1}{k + 1} \right) = 1 - \frac{1}{n + 1}\]

summation by parts gives

\[\sum_{k = 1}^{n}\frac{P_{k}}{k(k + 1)} = P_{n}s_{n} - \sum_{k = 1}^{n}{k^{2}\left| \alpha_{k} \right|^{2}s_{k - 1}}\]

\[= \sum_{k = 1}^{n}{k\left| \alpha_{k} \right|^{2} - \frac{1}{n + 1}\sum_{k = 1}^{n}{k^{2}\left| \alpha_{k} \right|^{2}}}\]

Thus

\[Q_{n} \leq (n + 1)\exp\left\{ \frac{1}{n + 1}\sum_{k = 1}^{n}{(n + 1 - k)\left( k\left| \alpha_{k} \right|^{2} - \frac{1}{k} \right)} \right\}\]

and the theorem follows.

 \end{proof}

~\\
We now consider again for any\(\ h \in S,\ \)where\(\ h\ \)is odd,
\(h(z) = \sqrt{f\left( z^{2} \right)},\ f \in S\). We define the
\emph{logarithmic coefficients} \(\ \gamma_{n}\ \)for a particular
function (for more details see \hyperlink{D1}{[D1]} p.151)

\[\log\frac{h\left( \sqrt{z} \right)}{\sqrt{z}} = \frac{1}{2}\log\left( \frac{f(z)}{z} \right) = \sum_{n = 1}^{\infty}{\gamma_{n}z^{n}}\]

With the above results we claim the following:

\uline{Conjecture 2.3.2 (Milin's conjecture)}

\fbox{\parbox{1.0\textwidth}{\ 
\[\sum_{m = 1}^{n}{\sum_{k = 1}^{m}\left( k\left| \gamma_{k} \right|^{2} - \frac{1}{k} \right)} \leq 0,\ n = 1,2,3\ldots\]

\emph{This conjecture implies Robertson's conjecture.}
}}

\begin{proof}

Exponentiating both sides of

\[\log\frac{h\left( \sqrt{z} \right)}{\sqrt{z}} = \frac{1}{2}\log\left( \frac{f(z)}{z} \right) = \sum_{n = 1}^{\infty}{\gamma_{n}z^{n}}\]

gives

\[\sum_{n = 0}^{\infty}{c_{2n + 1}z^{n} = \exp\left\{ \sum_{n = 1}^{\infty}{\gamma_{n}z^{n}} \right\}}\]

By Theorem 2.3.1 we have

\[\sum_{n = 0}^{\infty}{\left| c_{2n + 1} \right|^{2} \leq (n + 1)\exp\left\{ \sum_{m = 1}^{n}{\sum_{k = 1}^{m}\left( k\left| \gamma_{k} \right|^{2} - \frac{1}{k} \right)} \right\}} \leq (n + 1)\exp{(0) = n + 1}\]

which is Robertson's conjecture, and in turn proves Bieberbach's
conjecture.

 \end{proof}

\newpage 

\hypertarget{loewners-theory}{%
\section{3 Loewner's theory}\label{loewners-theory}}

Charles Loewner introduced the first non-elementary method in proving
Bieberbach's conjecture; this chapter will only be introducing the
theory partially, in order to prepare for Weinstein's proof of the
conjecture. Loewner's theory itself has many advanced results and
applications: readers are encouraged to turn to {[}D1{]} chapter 3,
{[}P1{]} Chapter 6, {[}S1{]} Chapter 16 and {[}K1{]} etc. for further
readings.

\subsection{3.1 Caratheodory's kernel convergence theorem
}
One fundamental concept in geometric function theory is \emph{kernel
convergence}. We first state the definition of this notion, and then
show that it is well-defined:

\uline {Definition 3.1.1 (\textit{Kernel Convergence})}:

Consider a sequence of domains\(\ D_{n},\ \)all containing the
origin(\(z = 0\)).

Case I:

Suppose\(\ \left\{ 0 \right\}\ \)is an interior point for
all\(\ D_{n}\), that is, if a disc\(\ |z| < r\ \)is contained in all the
domains. Then the \emph{kernel} is defined to be the largest domain,
call it\(\ D,\ \)such that every compact subset of\(\ D\ \)is contained
in all but a finite number of\(\ D_{n}\)s. That is, for every compact
(closed) subset of\(\ D,\ \)we can find a number\(\ N_{0}\ \)such that
the compact subset is contained in all\(\ D_{n}\ \)for\(\ n > N_{0}.\ \)

Case II:

Suppose\(\ \left\{ 0 \right\}\ \)is not an interior point for
all\(\ D_{n},\ \)that is, if the domains behave in such a way that no
matter how large \(n = n_{0}\ \)is chosen to be, there is no
disc\(\ |z| < r\ \)contained in all domains with\(\ n > n_{0}.\ \)Then
we define the kernel to be the origin itself.

In either case, we denote convergence to the
kernel\(\ D\ \)by\(\ D_{n} \rightarrow D.\)

Note that for a infinite sequence of strictly ``shrinking'' domains, it
is indeed possible for\(\ \left\{ 0 \right\}\ \)\emph{not} to be an
interior point: there is no open disc with radius small enough around it
such that the disc is contained in all domains as we can ``shrink'' it
even smaller!

Since this part was missed by most materials (for example {[}D1{]} and
{[}P1{]}), we show that \emph{kernels} are well-defined:

\begin{proof}

For Case I: The kernel is defined to be the largest domain with the
property: every compact subset of\(\ D\ \)is contained in all but a
finite number of\(\ D_{n}\)s.

Let\(\ P_{n}\ \)be a domain (not necessarily the largest) with this
property (that every compact subset of\(\ P_{n}\ \)is contained in all
but a finite number of\(\ D_{n}\)s, and
\(\left\{ 0 \right\} \subset P_{n}\)). If we could show that
\(\bigcup_{n = 1}^{\infty}{P_{n}\ }\)also has this property, then there
is indeed a \emph{largest} domain with this property because no larger
such domain exists, and hence proves that \emph{kernels} are
well-defined in this case.

Let\(\ U = \bigcup_{n = 1}^{\infty}{P_{n}\ }.\) Since
\({\left\{ 0 \right\} \subset P}_{n}\) for
all\(\ P_{n},\ \left\{ 0 \right\} \in U.\ \)For
each\(\ z_{0} \in U,\ \)there is some\(\ P_{n}\ \)such
that\(\ z_{0} \in P_{n},\ \)by definition of union. By the fact
that\(\ P_{n}\ \)is open, some open disc with positive radius is in it,
that is, \(\exists B\left( z_{0},r \right) \subset P_{n};\ \)by the
property of\(\ P_{n},\ \overline{B\left( z_{0},r \right)}\) is contained
in all but a finite number of\(\ D_{n}\)s, because a closed disc is
compact, and every compact subset of\(\ P_{n}\ \)satisfies this
property.

Now if\(\ K\ \)is a compact subset of\(\ U,\ \)then centred at every
point of\(\ K\ \)there is an open disc contained in
some\(\ P_{n},\ \)hence there exists an open cover of\(\ K.\ \)By the
Heine-Borel theorem, there is therefore a finite subcover
of\(\ K,\ \)which we could enumerate as\(\ B_{1},\ B_{2},\ B_{3}\ldots\)
Now for each of these open discs, its closure\(\ \overline{B_{n}}\) is
contained in all but a finite number of\(\ D_{n}\)s, as we proved above.
Therefore the entire of\(\ K \subset D_{n}\ \)for all but a finite
number of \(D_{n}\)s, because it is completely covered by these closed
discs, each contained in all but a finite number of\(\ D_{n}\)s. This
shows that \emph{kernels} are indeed well-defined for Case I.

For Case II: the kernel is defined to be a single point and does not
need the ``largest domain'' notion. This is clearly well-defined, and
hence completes the proof.

\end{proof}

We now introduce some lemmata to enhance our intuition for kernel
convergence.

\hypertarget{Ker}{\uline{Lemma 3.1.2 (Kernel of Increasing domains)}}

\fbox{\parbox{1.0\textwidth}{\ 
\emph{If}\(\ D_{n}\ \)\emph{is an increasing sequence of domains
containing}\(\ \left\{ 0 \right\},\ \)\emph{that is,}
\(D_{k} \subset D_{k + 1}\) \emph{for all}\(\ k,\ \)\emph{then the
kernel of this sequence of domains is its union.}
}}

\begin{proof}

By monotonicity, any subsequence of\(\ D_{n}\ \)converges to the same
limit, so it is sufficient to consider the sequence\(\ D_{n}\ \)itself.

We prove that the
kernel\(\ D \subseteq \bigcup_{n = 1}^{\infty}D_{n},\ \)and
also\(\ \bigcup_{n = 1}^{\infty}D_{n} \subseteq D,\ \)which
compels\(\ D\) to be equal to \(\bigcup_{n = 1}^{\infty}D_{n}.\)

First fix a point\(\ z_{0} \in D.\ \)By openness there is a closed disc
around this point that is fully contained in\(\ D.\ \)By definition of a
kernel, this disc (which is a compact subset) is contained in all but
finitely many\(\ D_{n}\)s. This means that there exists a
number\(\ n_{0}\ \)s.t. the disc is contained in
all\(\ D_{n}\ \)for\(\ n \geq n_{0},\ \)so there is a
domain\(\ {D_{n}}_{0}\ \)that fully contains it and all domains
after\(\ D_{n_{0}}\ \)also contains it. But since\(\ D_{n}\ \)is
increasing, this is impossible unless this disc is contained in
\(\bigcup_{n = 1}^{\infty}D_{n},\ \)by definition of an union. Since the
choice of this point is arbitrary, this shows
that\(\ D \subseteq \bigcup_{n = 1}^{\infty}D_{n}.\ \)

On the other hand, every point in an increasing sequence of domains
satisfies the property \(\ z_{0} \in D_{n} \Rightarrow z_{0} \in D_{n + 1}\). This means that
any compact subset of\(\ \bigcup_{n = 1}^{\infty}D_{n}\ \)also has this
property, as every point in this compact subset satisfies this property.
This already shows that \(\bigcup_{n = 1}^{\infty}D_{n} \subseteq D\)
and the proof is complete.

 \end{proof}

\uline{Lemma 3.1.3 (Decreasing domains)}

\fbox{\parbox{1.0\textwidth}{\ 

\emph{If}\(\ D_{n}\ \)\emph{is an decreasing sequence of domains
containing}\(\ \left\{ 0 \right\},\ \)\emph{that is,}
\(D_{k + 1} \subset D_{k}\) \emph{for all}\(\ k,\ \)\emph{then:}

~\\

\emph{1. If}\(\ \left\{ 0 \right\}\ \)\emph{is an interior point of the
intersection, that is, there is an open disc
around}\(\ \left\{ 0 \right\}\ \)

\emph{such that the disc is contained
in}\(\ \bigcap_{n = 1}^{\infty}D_{n},\ \)\emph{then the kernel is the
interior of} \(\bigcap_{n = 1}^{\infty}D_{n}.\ \)

~\\

\emph{2. If not, then the kernel is}\(\ \left\{ 0 \right\}.\ \)

}}

For a proof, see \hyperlink{P1}{[P1]} Lemma 1.4 p.28-29.

Now that we are more familiar with kernels and showed that there is no
problem with defining them this way, we are ready to introduce the next
theorem. Loewner's theory is based on the following important theorem,
attributed to Constantin Carathéodory:

\uline{Theorem 3.1.4 (Caratheodory's kernel convergence theorem):}

\fbox{\parbox{1.0\textwidth}{\ 

\emph{Let}\(\ \left\{ D_{n} \right\}\ \)\emph{be a sequence of simply
connected domains
with}\(\ \left\{ 0 \right\} \in D_{n}\mathbb{\subset C,\ }n = 1,2,\ldots\ \)

~\\ 

\emph{Let}\(\ f_{n}\left( \mathbb{D} \right) = D_{n},\ \)\emph{and
satisfies}\(\ f_{n}(0) = 0\ \)\emph{and}\(\ f_{n}^{'}(0) > 0.\ \)\emph{Let}\(\ D\ \)\emph{be
the kernel of}\(\ \left\{ D_{n} \right\}.\ \)\emph{Then}
\(f_{n} \rightarrow f\ \)\emph{uniformly on each compact subset
of}\(\mathbb{\ D\ }\)\emph{if and only
if}\(\ D_{n} \rightarrow D\mathbb{\neq C.\ }\)

~\\ 

\emph{Furthermore,
if}\(\ D = \left\{ 0 \right\},\ f(z) \equiv 0.\ \)\emph{Otherwise,}
\(\ f\left( \mathbb{D} \right) = D.\)

}}

We spend the rest of this chapter proving the theorem in two parts. We
will denote the kernel by\(\ D,\ \)and\(\ \varphi_{n}\ \)represents the
inverse function of\(\ f_{n},\ \)which always exists (locally) by
univalence of\(\ f_{n}.\ \)The proof follows from {[}D1{]} and {[}P1{]},
with supplementary details not given by most materials.

\begin{proof}

We begin with the proof of the first part: that is, if
\(f_{n}(z) \rightarrow f(z)\ \)uniformly on compact subsets
of\(\mathbb{\ D,\ }\)then\(\ D_{n} \rightarrow D.\ \)

\textbf{Part One: If}
\(\mathbf{f}_{\mathbf{n}}\left( \mathbf{z} \right)\mathbf{\rightarrow f}\left( \mathbf{z} \right)\mathbf{\ }\)\textbf{uniformly
on compact subsets
of}\(\mathbb{\ D,\ }\)\(\mathbf{\ }\mathbf{D}_{\mathbf{n}}\mathbf{\rightarrow D.\ }\)

Suppose that \(f_{n}(z) \rightarrow f(z)\ \)uniformly on compact subsets
of\(\mathbb{\ D.\ }\)Then\(\ f\ \)is analytic in\(\mathbb{\ D,\ }\)

and by Hurwitz's theorem (See \hyperlink{Hurwitz}{Theorem 1.1.2})\(\ f\ \)is either constant
or univalent. We divide into two cases for discussion.

Case 1: Suppose\(\ f\ \)is constant. \(f(0) = 0\ \)and\(\ f\ \)is
constant implies that \(f(z) \equiv 0.\ \)Suppose for the sake of
contradiction that\(\ D \neq \left\{ 0 \right\}.\ \)Then there exists a
simply connected domain\(\ P\ \)containing\(\ \left\{ 0 \right\}\ \)such
that each compact subset of\(\ P\ \)lies in all\(\ D_{n}\ \)except for a
finite number of\(\ D_{n}\)s.\(\ \)This implies that there is an open
disc\(\ B(0,\rho)\ \), where\(\ \rho\ \)is a positive number, that is a
subset of\(\ P,\ \)contained in all\(\ D_{n}\ \)for
\(n > N_{0},\ \)where\(\ N_{0}\ \)is a sufficiently large number.
Consider the inverse functions\(\ \varphi_{n}.\ \)Such inverse functions
are defined on the range of\(\ f_{n},\ \)in this case we consider the
set\(\ B(0,\rho)\), which is a subset of the range of\(\ f_{n}.\ \)In
this disc, \(\varphi_{n}(0) = 0\ \)and since \(f_{n}(z)\ \)is defined
on\(\mathbb{\ D,\ }\left| \varphi_{n}(z) \right| < 1.\ \)

By the Schwarz lemma,
\(\left| \varphi^{'}(0) \right| \leq \frac{1}{\rho}.\ \)By the
well-known fact that derivatives of inverse functions are reciprocals
(at corresponding points), and
that\(\ f_{n}(0) = \varphi_{n}(0) = 0,\ \)we have that for sufficiently
large\(\ n,\ f_{n}^{'}(0) \geq \rho.\) But this means that
\(\ f_{n}(z) \nrightarrow 0.\ \)Hence\(\ D = \left\{ 0 \right\}.\ \)
Replacing\(\ f_{n}\ \)with an arbitrary subsequence shows that this
result applies to any subsequence of\(\ D_{n},\ \)hence the kernel is
indeed by definition\(\ \left\{ 0 \right\}.\ \)

Case 2:

Suppose\(\ f\ \)is not constant.
Then\(\ f\left( \mathbb{D} \right) = \Delta\ \)where\(\ \Delta\ \)is
some domain. We prove that\(\ \Delta\ \)coincides with\(\ D,\ \)by first
showing\(\ \Delta \subseteq D\ \)and then
also\(\ D \subseteq \Delta,\ \)which in turn shows
that\(\ \Delta = D.\ \)

Step I: \(\Delta \subseteq D\)

For any arbitrary compact subset of\(\ \Delta\ \)the inverse
function,\(\ f^{- 1}(z),\ \)is defined in it, because it is a subset of
the image of\(\ f(z).\ \)Choose an arbitrary compact subset, call
it\(\ E,\ \)and surround it with a rectifiable (that is, it is closed
and does not include infinity as a point) Jordan
curve\(\ \Gamma,\ \)where\(\ \Gamma \in \left\{ \Delta \right\}\backslash\text{\{}E\}\ .\ \)Now
let\(\ \gamma\ \)be the pre-image
of\(\ \Gamma\ \)under\(\ \ f(z),\ \)that
is,\(\ f(\gamma) = \Gamma.\ \)We prove that\(\ E \in D_{n}\ \)for
sufficiently large\(\ n.\ \)

Indeed, fixing a point\(\ w_{0} \in E,\ \)then by the fact that
\(E\ \)belongs to the image of the interior of \(\gamma\ \)(See
\hyperlink{Int}{Proposition 1.3.2}) we have
that\(\ \left| f(z) - w_{0} \right| \geq \delta\ \)for
all\(\ z \in \gamma;\ \)that is, \(f(\gamma)\ \)always keeps a positive
distance with the arbitrary point we fixed in\(\ E.\ \)Since we are
assuming\(\ f_{n} \rightarrow f\ \)uniformly, we have that there exists
some
number\(\ N\ \)s.t.\(\ n \geq N \Longrightarrow \left| f_{n}(z) - f(z) \right| < \delta,\ \)by
definition of uniform convergence for sequences of functions.

This means that for\(\ n \geq N,\ \)

\[\left| f_{n}(z) - f(z) \right| = \left| \left( f_{n}(z) - w_{0} \right) - \left( f(z) - w_{0} \right) \right| < \delta \leq \left| f(z) - w_{0} \right|\]

We compare the functions \(\left| f_{n}(z) - f(z) \right|\ \)and
\(\left| f(z) - w_{0} \right|.\ \)By Rouches' theorem, noticing that
inside\(\ \gamma\ \)it holds that
\(\left| f_{n}(z) - f(z) \right| \leq \left| f(z) - w_{0} \right|,\ \)we
can conclude that
\(\left\lbrack f(z) - w_{0} \right\rbrack\ \)and\(\ \left\lbrack f_{n}(z) - f(z) \right\rbrack + \left\lbrack f(z) - w_{0} \right\rbrack\ \)has
the same number of zeroes. But
\(\left\lbrack f_{n}(z) - f(z) \right\rbrack + \left\lbrack f(z) - w_{0} \right\rbrack = \left\lbrack f_{n}(z) - w_{0} \right\rbrack\),
hence
\(\left\lbrack f(z) - w_{0} \right\rbrack\ \)and\(\ \left\lbrack f_{n}(z) - w_{0} \right\rbrack\ \)has
the same number of zeroes. Now\(\ f(z) = w_{0}\ \)at exactly one point
by its univalence, so there is only one zero for
\(\left\lbrack f(z) - w_{0} \right\rbrack,\ \)and thereby one zero for
\(\left\lbrack f_{n}(z) - w_{0} \right\rbrack.\ \)That is, for any
arbitrary\(\ w_{0} \in E\ \)that we choose, it is
inside\(\ f_{n}\left( \mathbb{D} \right) = D_{n}\ \)for\(\ n \geq N,\ \)and
hence\(\ E\ \)is inside all\(\ D_{n}\)s for sufficiently
large\(\ n.\ \)By definition of \emph{kernels},
\(\Delta \subseteq D,\ \)because every compact subset of it belongs to
all\(\ D_{n}\)s for sufficiently large\(\ n.\)

Step II: \(D \subseteq \Delta\)

Clearly\(\ \left\{ 0 \right\} \in \Delta = f\left( \mathbb{D} \right)\ \)as\(\ f(0) = 0.\ \)Choose
an arbitrary point\(\ z_{0} \neq 0.\ \)Consider a closed disc\(\ B\)
around\(\ z_{0}\ \)that is fully contained in\(\ D.\) Since every
compact subset of\(\ D\ \)belongs to all but finitely many domains, this
disc belongs to all\(\ D_{n}\)s
for\(\ n > n_{0},\ \)where\(\ n_{0}\ \)is some finite number. Consider
the domain that is equivalent to the interior of\(\ B.\ \)The inverse
functions\(\ \varphi_{n}(z)\ \)are analytic in this disc, because it is
fully contained in its domain (which is the range of\(\ f_{n}(z)\)) and
satisfies\(\ \left| \varphi_{n}(z) \right| < 1\ \)(as its range is the
domain of\(\ f_{n}(z)\ \)which is the unit disc), for
\(n > n_{0}.\ \)Since this sequence of functions are analytic and
bounded, it is a normal family by Montel's theorem (See \hyperlink{Montel}{Theorem 1.1.5}).
This means that there exists a
subsequence\(\ \ {\varphi_{n}}_{k}\ \)that converges locally uniformly
in the interior of\(\ B,\ \)given that \(n_{k} > n_{0}.\) Since it is
uniformly convergent and satisfies\(\ \varphi_{n}(0) = 0,\ \)the limit
function\(\ \varphi(z)\ \) (that is, the function
satisfying\(\ \varphi_{n_{k}}(z) \rightarrow \varphi(z)\))
satisfies\(\ \varphi(0) = 0\ \)and\(\ \left| \varphi(z) \right|1.\ \)Since\(\ \varphi_{n_{k}}(z)\ \) converges locally uniformly around\(\ w_{0},\ f_{n_{k}}(z)\ \)converges
locally uniformly around \(\ \varphi\left( z_{0} \right).\ \)(We make the remark that connected
sets are always mapped to connected sets under analytic univalent
functions, See \hyperlink{Int}{Proposition 1.3.2}).
Now\(\ \varphi\left( z_{0} \right)\ \)is some point in the unit disc,
and it holds
that\(\ \varphi_{n_{k}}\left( z_{0} \right) \rightarrow \varphi\left( z_{0} \right)\),
so it is easy to see that\(\ z_{0} \in \Delta:\ \)indeed,
suppose\(\ \varphi\left( z_{0} \right) = z_{1},\ \)and we know that
\(f_{n_{k}}\left( z_{1} \right) \rightarrow f\left( z_{1} \right)\)
around\(\ z_{1},\ \)then based on the fact
that\(\ f_{n_{k}}\left( z_{1} \right) \in \Delta\forall n_{k}\ \)and
uniform convergence, \(f\left( z_{1} \right) \in \Delta.\ \)

Based on the fact that\(\ D \subseteq \Delta\ \)and
also\(\ \Delta \subseteq D\ \)we conclude
that\(\ f\left( \mathbb{D} \right) = D\)

\textbf{Part Two: If}
\(\mathbf{D}_{\mathbf{n}}\mathbf{\rightarrow D}\)\textbf{, then}
\(\mathbf{f}_{\mathbf{n}}\left( \mathbf{z} \right)\mathbf{\rightarrow f}\left( \mathbf{z} \right)\mathbf{\ }\)\textbf{uniformly
on compact subsets of}\(\mathbb{\ D}\)\textbf{.}

Now suppose\(\ D_{n} \rightarrow D.\ \)By the Koebe one-quarter theorem
(See \hyperlink{Koebe}{Theorem 1.3.7}), it is easy to see that\(\ f^{'}(0)\ \)is bounded,
or else the disc\(\ |z| < \frac{1}{4}f_{n}^{'}(0)\ \)belongs to
all\(\ D_{n}\)s which in turn sets the limit domain to the entire
of\(\mathbb{\ C,\ }\)contradicting the definition of a kernel. 

Moreover,
by the growth theorem (See \hyperlink{Growth}{Theorem 1.4.3}) it holds that

\[\left| f_{n}(z) \right| \leq f^{'}(0)\frac{|z|}{\left( 1 - |z| \right)^{2}}\]

and hence it follows that \(f_{n}(z)\ \)is bounded locally uniformly. By
Montel's theorem (See \hyperlink{Montel}{Theorem 1.1.5}), this means that\(\ f_{n}(z)\ \)is
a normal family, and hence there exists a subsequence that converges for
every sequence we choose.

Suppose for the sake of contradiction that\(\ f_{n}(z)\ \)does not
converge locally uniformly. Then there exists two
sub-sequences\(\ f_{n_{k}}\ \)and\(\ {f_{n}}_{l}\ \)converging to
different limits\(\ k(z)\ \)and\(\ l(z)\), by the fact
that\(\ f_{n}\ \)constitutes a normal family.

But from \textbf{Part One} this means that the limit functions for these
two sub-sequences maps

\(\mathbb{\ D}\) to different domains: indeed, in view of the Riemann
mapping theorem, the fact that \(k(0) = l(0) = 0\ \)and\(\ k^{'}(0) \geq 0,\ l^{'}(0) \geq 0\ \)guarantees
that if they map\(\mathbb{\ D\ }\)to different domains, then they are
identical because the mapping is \emph{unique} for any domain
in\(\mathbb{\ C}\). This contradicts to the assumption
that\(\ D_{n} \rightarrow D\ \)by the definition of kernels.

With the above two \textbf{Parts} combined, the proof of this theorem is
now complete. 

 \end{proof}

\subsection{3.2 Loewner's theory and Further Lemmas}

We shall first introduce the theory of subordinating chains. To do so we
need the following special classes of functions, each easy to define:

\textbf{The Caratheodory class of functions:}

Functions that are holomorphic
in\(\mathbb{\ D,\ }\)sending\(\ 0\ \)to\(\ 1\ \)and has positive real
part. We denote this class with\(\ p:\ \ \ \ p(0) = 1\) and
\(Re(p) > 0.\ \)

\textbf{The Schwarz functions:}

Functions that are holomorphic in\(\mathbb{\ D,\ }\)sending zero to zero
and has modulus smaller than 1. We denote this class
with\(\ \varphi:\ \ \ \ \ \varphi(0) = 0\ \)and\(\ \left| \varphi(z) \right| < 1.\ \)

\uline{Definition 3.2.1 (Subordination):}

\fbox{\parbox{1.0\textwidth}{\ 
We say\(\ f\ \)is \emph{subordinate}
to\(\ g\ \)if\(\ f(z) = g\left( \varphi(z) \right)\ \)for some Schwarz
function\(\ \varphi.\ \)

}}

In this case we denote this by

\[f(z) \prec g(z)\]

Notice that a function belongs to the Caratheodory class if and only if
it can be written in the form

\[p(z) = \frac{1 + \varphi}{1 - \varphi}\]

For some Schwarz function\(\ \varphi\). To see this we complete the following exercise:

\newpage

\uline{Theorem 3.2.2 (Subordination Principle)}

\fbox{\parbox{1.0\textwidth}{\ 

\emph{If}\(\ f\ \)\emph{is univalent in}\(\mathbb{\ D\ }\)\emph{and
if}\(\ g\ \)\emph{is a function analytic
in}\(\mathbb{\ D\ }\)\emph{with}\(\ g(0) = f(0)\ \)\emph{and}\(\ \)

\(g\left( \mathbb{D} \right) \subset f\left( \mathbb{D} \right),\ \)\emph{then}\(\ \left| g^{'}(0) \right| \leq \left| f^{'}(0) \right|.\ \ \)

~\\

\emph{Also it holds
that}\(\ g\left(\mathbb{D}_{r} \right)\subset f\left(\mathbb{D}_{r} \right),\ \)\emph{where}\(\ r < 1.\ \)\emph{We hence say that}\(\ g\ \)\emph{subordinates to}\(\ f.\ \)
}}

\begin{proof}

WLOG we normalize \(g(0) = f(0) = 0.\) The theorem then follows as we
consider transformations on the functions\(\ g\ \)and\(\ f.\ \) Consider
the function\(\ {P(z) = f}^{- 1}\left( g(z) \right).\ \)By univalence
of\(\ f\ \)we know that\(\ f^{- 1}\ \)exists and since
\(g\left( \mathbb{D} \right) \in f\left( \mathbb{D} \right),\ \)the
domain of\(\ P\ \)is well defined. Applying the Schwarz Lemma
on\(\ P\ \)(since\(\ P\ \)satisfies\(\ {P(0) = f}^{- 1}\left( g(0) \right) = f^{- 1}(0) = 0\ \)and\(\ \left| P(z) \right| < 1\ \)as
\(P(z)\mathbb{\in D)}\) we have
that\(\ \left| P^{'}(0) \right| \leq 1\ \)and\(\ \left| P(z) \right| \leq |z|\).
It follows
that\(\ \left| P^{'}(0) \right| = \left| \frac{g^{'}(0)}{f^{'}(0)} \right| \leq 1\)
by the chain rule
and\(\ \left| f^{- 1}\left( g\left( D(0,r) \right) \right) \right| \leq \left| D(0,r) \right|.\ \)The
theorem then follows.

It only remains to notice that

\[p(z) \prec h(z) = \frac{1 + z}{1 - z}\]

As \(p(0) = 1 = h(0)\).
Also\(\ p\left( \mathbb{D} \right) \subset h\left( \mathbb{D} \right)\ \)because
the later maps the unit disc onto the right half plane and the former
has positive real part. Hence it follows that

\[p(z) \prec \frac{1 + z}{1 - z} \Longrightarrow p(z) = \frac{1 + \varphi}{1 - \varphi}\]

For some Schwarz function\(\ \varphi\), which is in accord with the
definition of Subordination.

\end{proof}

We now make an introduction to Loewner's Theory

\uline{Definition 3.2.3 (\textit{Slit Mappings}):}

\fbox{\parbox{1.0\textwidth}{\ 
We say\(\ f\ \)is \emph{slit} if it maps a domain conformally onto the
entire complex plane except a set of Jordan arcs. A \emph{single-slit
mapping} is a slit mapping whose range is the complement of a
\emph{single} Jordan arc.

}}

Loewner's theory concerns approximating functions in the
class\(\ S\ \)uniformly by single-slit mappings, and monitoring its
dynamics with differential representations, for example the well-known
Loewner's differential equations. To do so, we first need to proof that
indeed,

\uline{Theorem 3.2.4}

\fbox{\parbox{1.0\textwidth}{\ 

\emph{For each}\(\ f \in S\ \)\emph{there corresponds a sequence of
single-slit mappings}\(\ f_{n} \in S\ \)\emph{such that} \(f_{n} \rightarrow f\ \)\emph{uniformly on each compact subset
of}\(\mathbb{\ D.\ }\)

}}

\begin{proof}

Consider univalent maps\(\ g_{n}(z)\) that maps\(\mathbb{\ D\ }\)to the
interior of an analytic Jordan curve\(\ \Gamma.\ \) In view of the
Riemann Mapping theorem, we assume that all\(\ g_{n}\ \)satisfies the
conditions \(\ g(0) = 0,\ g^{'}(0) > 0.\ \) Observe that there exists a
sequence of\(\ g_{n}\) that converges uniformly to
each\(\ f \in S,\ \)for
\(\frac{f_{n}\left( r_{n}z \right)}{r_{n}},\ 0 < r < 1,r = 1 - 1/n\ \)
gives an example of such a sequence.

Consider an arbitrary function\(\ f \in S\ \)that
maps\(\mathbb{\ D\ }\)to the interior of a Jordan curve\(\ \Gamma.\ \)

Let\(\ \Gamma_{n}\ \)be another Jordan arc connecting the point infinity
to a point\(\ w_{0} \in \Gamma,\ \)then follows the exact same path
as\(\ \Gamma\ \)starting from\(\ w_{0},\ \)until it stops at some
point\(\ w_{n}.\ \)Let\(\ D_{n}\ \)be the complement
of\(\ \Gamma_{n}\ \)(that is,
\(D_{n} = \mathbb{C}^{\infty}\backslash\Gamma_{n}\)) and
let\(\ g_{n}\ \)map\(\mathbb{\ D\ }\)conformally onto\(\ D_{n},\ \) with
the conditions \(g_{n}(0) = 0,\ g_{n}^{'}(0) > 0.\ \)Such
functions\(\ g_{n}\ \)exists in view of the Riemann Mapping theorem. Let
the endpoints,\(\ w_{n},\) of\(\ \Gamma_{n}\) be chosen such
that\(\ \Gamma_{n} \subset \Gamma_{n + 1}\ \)and\(\ w_{n} \rightarrow w_{0}.\ \)

We make the remark that as\(\ t_{n} \rightarrow T < \infty,\ \)the
kernel of the sequence of domains\(\ D_{t_{n}}\ \)is\(\ D_{T}.\ \)This
is because the domains are strictly increasing, and increasing sequences
of domains converges to its union. (See \hyperlink{Ker}{Lemma 3.1.2}). That is, we have
that the kernel of the sequence of domains\(\ D_{n}\ \)is\(\ D.\ \)

Therefore by Caratheodory's convergence
theorem,\(\ g_{n} \rightarrow f\ \)uniformly on compact subsets
of\(\mathbb{\ D.\ }\)We know that for uniformly convergent sequences of
functions, \(g_{n}^{'}(0) \rightarrow f^{'}(0) = 1\ \)by Cauchy's
integral formula, hence

\[h_{n} = \frac{g_{n}(z)}{g_{n}^{'}(0)} \in S\]

are single-slit mappings that converges uniformly to\(\ f\ \)on compact
subsets of\(\mathbb{\ D.\ }\)

\end{proof}

With the above result, we now produce a system to aid us in representing
the single-slit mappings that goes uniformly to any\(\ f \in S.\ \)We
follow the footpath of {[}D1{]} p.82.

Let\(\ f \in S\ \)be a single-slit mapping from\(\mathbb{\ D\ }\)onto a
domain\(\ D\ \)which is the complement of a Jordan
arc\(\ \Gamma\ \)extending from a point\(\ w_{0}\mathbb{\in C\ }\)to the
point infinity. We parametrize\(\ \Gamma\ \)by a
function\(\ \psi(t),\ 0 \leq t < T\) (we cannot
conclude\(\ T = \infty\ \)just yet, but we may do so later) such
that\(\ \psi(0) = w_{0}\ \)and\(\ s \neq t \Longrightarrow \psi(s) \neq \psi(t).\ \)

Let\(\ \Gamma_{t}\ \)be the portion
of\(\ \Gamma\ \)from\(\ \psi(t)\ \)to\(\ \infty,\ \)and
let\(\ D_{t} = \mathbb{C}^{\infty}\backslash\Gamma_{t}.\)
Then\(\ D_{0} = D\ \)and\(\ D_{s} \subset D_{t}\ \)for\(s < t.\ \)To visualize this better, perhaps think of\(\ \Gamma_{t}\ \)as a fire fuse: as\(\ t \rightarrow T,\ \)the acr\(\ \Gamma_{t}\ \)``burns''
from one end (starting from\(\ w_{0}\)) all the way to the point
\(\infty\), leaving the remaining part of the plane (\(D_{t}\)) larger.

Let

\[g(z,t) = \beta(t)\left\{ z + b_{2}(t)z^{2} + \ldots \right\}\]

be the series expansion of the conformal mapping
of\(\mathbb{\ D\ }\)onto\(\ D_{t}\ \)and, in view of the Riemann Mapping
theorem, we may
let\(\ g(0,t) = 0\ \)and\(\ g^{'}(z,t) = \beta(t) > 0.\ \)

Now fix\(\ z\ \)and consider the function\(\ g(z,t)\ \)as a function
of\(\ t.\ \)Then it is easy to see that by Caratheodory's theorem,
\(g_{t_{n}} \rightarrow g_{t}\ \) locally uniformly
as\(\ {\ t}_{n} \rightarrow t.\ \)

In particular, \(\beta\left( t_{n} \right) \rightarrow \beta(t),\ \)
locally uniformly for any sequence of\(\ t_{n}\ \)arbitrarily
chosen\(.\ \)By the sequential definition of continuity, it follows
that\(\ \beta(t)\ \)is continuous in\(\ t\). This result need not rely
on Cauchy's integral formula for coefficients ({[}P1{]} p.156). However,
to see that the rest of the coefficients are also continuous
coefficients of\(\ t,\ \)notice that by Cauchy's formula for
coefficients,

\[a_{n}(t) = \frac{1}{2\pi i}\int_{C}^{}{\frac{g_{t}(\xi)}{\xi^{n + 1}}d\xi}\]

and continuity of\(\ a_{n}(t)\ \)follows from the differentiability
of\(\ g.\)

It is clear that\(\ \beta(0) = 1,\ \)as\(\ g(z,0) = f(z).\ \)This can be
seen by considering the fact that \(D_{0} = D\), so\(\ g(z,0)\ \)simply
maps\(\mathbb{\ D\ }\)to\(\ D.\ \)In view of the definition of
subordination (combined with the geometric view of the images), it is
easy to see that\(\ g(z,s) \prec g(z,t)\ \)

for\(\ s \leq t.\ \)So by the Subordination Principle (Theorem 3.2.2),
we know that

\[\left| g^{'}(0,s) \right| \leq \left| g^{'}(0,t) \right| \Longrightarrow \left| \beta(s) \right| \leq \left| \beta(t) \right|\]

and since\(\ \beta(t) > 0\ \)is a real function in\(\ t,\ \)not a
complex function, we conclude that it is increasing. It is easy to see
that it is \emph{strictly} increasing: otherwise by the uniqueness
guaranteed by Riemann Mapping theorem there exists \(D_{s} = D_{t}\),
\(s \neq t,\ \)contradicting to our definition of the
domains\(\ D_{t}\).

We may reparametrize \(\Gamma\ \)such
that\(\ \beta(t) = e^{t},\ 0 \leq t < T\): this is a standard approach
in analytic geometry (See \hyperlink{A2}{[A2]} p.10 chapter 1.3). We hence have
chosen a parametrization of the omitted arc\(\ \Gamma\ \)such that
(renaming\(\ g\ \)as\(\ f\) for distinction)

\[f(z,t) = e^{t}\left\{ z + \sum_{n = 2}^{\infty}{b_{n}(t)z^{n}} \right\}\]

where\(\ 0 \leq t < T.\ \)

We now show that actually it is necessary for \(T\ \)to be\(\ \infty.\)

We fix\(\ M \in \mathbb{R}^{+}\ \), and notice that there must be a
circle\(\ |w| = M\ \)such that the omitted arc \(\Gamma_{t}\ \)lies entirely outside this circle for sufficiently
large\(\ t\ \)(indeed, again with the fire fuse visualization, the fuse
must burn continuously to\(\ \infty\ \)and will eventually have no
remaining part in any circle we draw on the plane). From the maximum
modulus principle,

\[\left| \frac{z}{g(z,t)} \right| \leq \frac{1}{M},\ |z| < 1\]

Rearrange to get, since\(\ g(0,t) = 0,\ \)

\[M \leq \left| \frac{g(z,t)}{z} \right| = \left| \frac{g(z,t) - g(0,t)}{z - 0} \right|\]

By definition of derivatives,
taking\(\underset{z \rightarrow 0}{\ lim}\) on both sides we arrive at
\(M \leq \left| g^{'}(0,t) \right| = e^{t}\ \)for all\(\ t\ \)close
enough to\(\ T.\ \)Since\(\ M\ \)is arbitrary, this shows that
\(e^{t} \rightarrow \infty\ \)as\(\ t \rightarrow T,\ \ \)so it must
hold that\(\ T = \infty.\ \) 

With the above result, we may derive a differential representation of
the functions\(\ f_{t}.\ \) 

Although this representation was used by de Branges in his proof to
Bieberbach's conjecture, Weinstein's thesis did not make use of it, so
we shall not discuss Loewner's differential equations in details.

\newpage 

\uline{Theorem 3.2.5(Loewner's Differential Equations)}

\emph{For each}\(\ f \in S\ \)\emph{be a single slit mapping with
omitted arc}\(\ \Gamma.\ \)\emph{Then there exists a continuous
function}\(\ \kappa(t):\lbrack 0,\infty) \rightarrow \partial\mathbb{D\ }\)\emph{and
a family}\(\ f_{t}(z) = g^{- 1}\left( f(z),t \right)\ \)\emph{such that}

\[\frac{\partial f_{t}}{\partial t} = - f_{t}\frac{1 + \kappa(t)f_{t}}{1 - \kappa(t)f_{t}}\]

\emph{and}\(\ \lim_{t \rightarrow \infty}{e^{t}f_{t}(z) = f(z)}\)\emph{.}

We still outlined this statement, since it is a significant result, and
may interest readers to study its applications: for further
investigations of Loewner's theory and a proof to the above statement,
readers are encouraged to consult \hyperlink{D1}{[D1]}p.83, \hyperlink{S1}{[S1]} p.241 Theorem
16.2.2.

We make the remark that we have already shown that for
any\(\ f \in S\ \)there exists a family\(\ f_{t}\ \)of single-slit
mappings such that:

\begin{enumerate}
\def\labelenumi{\arabic{enumi}.}
\item
  \(\ f_{0} = f\)
\item
  \(f_{t}(z) = e^{t}z + a_{2}(t)z^{2} + \ldots\)
\end{enumerate}

Notice also that, considering the geometric interpretation of this
family,

\[f(z,s) \prec f(z,t)\]

Which implies by definition of subordination that

\[f(z,s) = f\left( \varphi(z,s,t),t \right)\]

We next prove the following:

\[Re\left\{ \frac{\frac{\partial f_{t}(z)}{\partial t}}{z\frac{\partial f_{t}(z)}{\partial z}} \right\} > 0\]

\newpage

To do so we first prove the following lemmata:

\uline{Lemma 3.2.6}

\[\left| f(z,t) - f(z,s) \right| \leq \frac{8|z|}{\left( 1 - |z| \right)^{4}}\left( e^{t} - e^{s} \right)\]

\uline{Lemma 3.2.7}

\[\left| \varphi(z,t,u) - \varphi(z,s,u) \right| \leq \frac{2|z|}{\left( 1 - |z| \right)^{2}}\left( 1 - e^{s - t} \right)\]

\begin{proof}

Consider first the function defined by

\[p(z,s,t) = \frac{1 + e^{s - t}}{1 - e^{s - t}}\frac{1 - z^{- 1}\varphi(z,s,t)}{1 + z^{- 1}\varphi(z,s,t)}\]

Now\(\ z^{- 1}\varphi(z,s,t)\ \)has Taylor
expansion\(\ e^{s - t} + \ldots\ \)because\(\ f\ \)has Taylor expansion

\(e^{t}z + \ldots\)and \(f(z,s) = f\left( \varphi(z,s,t),t \right)\). So

\[p(0,s,t) = \frac{1 + e^{s - t}}{1 - e^{s - t}}\frac{1 - e^{s - t}}{1 + e^{s - t}} = 1\]

It is easy to see that the function\(\ p(z,s,t)\) has positive real
part. (The proof was omitted by most materials)

Indeed,
suppose\(\ w_{s,t} = \varphi(z,s,t) = \varphi(p + qi,s,t) = x + yi\).
Since \(\frac{1 + e^{s - t}}{1 - e^{s - t}}\ \)is a positive number, we
consider only

\[\frac{1 - z^{- 1}\varphi(z,s,t)}{1 + z^{- 1}\varphi(z,s,t)} = \frac{z - \varphi}{z + \varphi}\]

Now

\[\frac{z - \varphi}{z + \varphi} = \frac{p + qi - x - yi}{p + qi + x + yi}\]

\[= \frac{(p - x) + (q - y)i}{(p + x) + (q + y)i}\]

\[= \frac{\left( p^{2} + q^{2} \right) - \left( x^{2} + y^{2} \right) + 2(xq - py)i}{(p + x)^{2} + (q + y)^{2}}\]

\[= \frac{\left| z^{2} \right| - \left| w^{2} \right| + 2(xq - py)i}{(p + x)^{2} + (q + y)^{2}}\]

Which in turn shows that

\[Re\left( \frac{z - \varphi}{z + \varphi} \right) = \frac{\left| z^{2} \right| - \left| w^{2} \right|}{(p + x)^{2} + (q + y)^{2}} \geq 0\]

as \(\left| z^{2} \right| - \left| w^{2} \right| \geq 0\ \)by
subordination.

Hence\(\ p(z,s,t)\ \)belongs to the Caratheodory class. It follows that

\[\ \ p(z,s,t) = \frac{1 + \varphi}{1 - \varphi} \Longrightarrow \left| p(z,s,t) \right| \leq \frac{1 + |z|}{1 - |z|}\]

\[\Longrightarrow \frac{1 + e^{s - t}}{1 - e^{s - t}}\frac{\left| z - \varphi(z,s,t) \right|}{\left| z + \varphi(z,s,t) \right|} \leq \frac{1 + |z|}{1 - |z|}\]

\[\Longrightarrow \left( 1 + e^{s - t} \right)\left| z - \varphi(z,s,t) \right| \leq \frac{1 + |z|}{1 - |z|}\left| z + \varphi(z,s,t) \right|\left( 1 - e^{s - t} \right)\]

\[\Longrightarrow \left| z - \varphi(z,s,t) \right| \leq \left( 1 + e^{s - t} \right)\left| z - \varphi(z,s,t) \right| \leq \frac{1 + |z|}{1 - |z|}\left| z + \varphi(z,s,t) \right|\left( 1 - e^{s - t} \right)\]

\[\Longrightarrow \left| z - \varphi(z,s,t) \right| \leq \frac{1 + |z|}{1 - |z|}|z + z|\left( 1 - e^{s - t} \right)\]

as\(\ |\varphi| \leq |z|\)

\[\Longrightarrow \left| z - \varphi(z,s,t) \right| \leq 2|z|\frac{1 + |z|}{1 - |z|}\left( 1 - e^{s - t} \right)\]

By the Distortion Theorem (See \hyperlink{Dist}{Theorem 1.4.2})

\[\left| f^{'}(\xi,t) \right| \leq \frac{2e^{t}}{\left( 1 - |\xi| \right)^{3}}\]

So it follows from the subordinating property that

\[\left| f(z,t) - f(z,s) \right| = \left| \int_{\varphi(z,s,t)}^{z}{f^{'}(\xi,t)d\xi} \right| \leq \left| z - \varphi(z,s,t) \right|\frac{2e^{t}}{\left( 1 - |z| \right)^{3}}\]

\[\leq \frac{2e^{t}}{\left( 1 - |z| \right)^{3}}2|z|\frac{1 + |z|}{1 - |z|}\left( 1 - e^{s - t} \right) = \frac{4|z|\left( e^{t} - e^{s} \right)\left( 1 + |z| \right)}{\left( 1 - |z| \right)^{4}} \leq \frac{8|z|\left( e^{t} - e^{s} \right)}{\left( 1 - |z| \right)^{4}}\]

which is Lemma 3.2.6.

Lemma 3.2.7 can be proven in the exact same manner, replacing the
integrand from \(f^{'}(\xi,t)\ \)to\(\ \varphi^{'}(\xi,t,u)\ \)and noticing
that\(\ |\varphi'| \leq \left( 1 - |z|^{2} \right)^{- 1}.\ \)

It only remains to notice that, by rearranging,

\[\frac{f(z,t) - f(z,s)}{t - s} = \left( \frac{e^{t - s} - 1}{e^{t - s} + 1} \cdot \frac{1 + e^{s - t}}{1 - e^{s - t}} \right)\left( \frac{z + \varphi}{z - \varphi} \cdot \frac{z - \varphi}{z + \varphi} \right)\frac{f(z,t) - f(\varphi,t)}{t - s}\]

\[= \left( \frac{1 + e^{s - t}}{1 - e^{s - t}} \cdot \frac{z - \varphi}{z + \varphi} \right)\frac{f(z,t) - f(\varphi,t)}{z - \varphi}\frac{z + \varphi}{e^{t - s} + 1}\frac{e^{t - s} - 1}{t - s}\]

\[= p(z,s,t)\frac{f(z,t) - f(\varphi,t)}{z - \varphi}\frac{z + \varphi}{e^{t - s} + 1}\frac{e^{t - s} - 1}{t - s}\]

By Lemma 3.2.6 it holds
that\(\ f(z,t) \rightarrow f(z,s)\ \)as\(\ t \rightarrow s\ \)locally
uniformly in\(\ D\ \)and hence \(f^{'}(z,t) \rightarrow f^{'}(z,s)\ \)(follows from estimation lemma).

Now

\[\frac{\partial f\left( z + \lambda\lbrack\varphi - z\rbrack,t \right)}{\partial\lambda} = (\varphi - z)f^{'}\left( z + \lambda\lbrack\varphi - z\rbrack,t \right)\]

\[\Longrightarrow f\left( z + \lambda\lbrack\varphi - z\rbrack,t \right) = \int_{}^{}{(\varphi - z)f^{'}\left( z + \lambda\lbrack\varphi - z\rbrack,t \right)}d\lambda\]

\[\Longrightarrow (\varphi - z)\int_{1}^{0}{f^{'}\left( z + \lambda\lbrack\varphi - z\rbrack,t \right)}d\lambda = f(z,t) - f(\varphi,t)\]

\[\Longrightarrow \int_{0}^{1}{f^{'}\left( z + \lambda\lbrack\varphi - z\rbrack,t \right)}d\lambda = \frac{f(z,t) - f(\varphi,t)}{z - \varphi}\]

But by Lemma 3.2.7, as\(\ t \rightarrow s,\ \varphi \rightarrow z.\ \)So

\[\frac{f(z,t) - f(\varphi,t)}{z - \varphi} = \int_{0}^{1}{f^{'}\left( z + \lambda\lbrack\varphi - z\rbrack,t \right)}d\lambda \rightarrow \int_{0}^{1}{f^{'}(z,s)}d\lambda = f^{'}(z,s)\]

It follows that letting \(t \rightarrow s\),

\[\frac{f(z,t) - f(z,s)}{t - s} \rightarrow \frac{\partial f(z,s)}{\partial s}\]

and hence

\[p(z,s,t)\frac{f(z,t) - f(\varphi,t)}{z - \varphi}\frac{z + \varphi}{e^{t - s} + 1}\frac{e^{t - s} - 1}{t - s} \rightarrow p(z,s,t)f^{'}(z,s)\frac{2z}{2} = zp(z,s,t)f^{'}(z,s)\]

So (this is referred to as \emph{Loenwer's PDE})

\[p(z,s,t) = \frac{\frac{\partial f_{s}(z)}{\partial s}}{z\frac{\partial f_{s}(z)}{\partial z}}\]

And the fact that

\[Re\left\{ \frac{\frac{\partial f_{t}(z)}{\partial t}}{z\frac{\partial f_{t}(z)}{\partial z}} \right\} > 0\]

Follows from the positivity of the real part of\(\ p(z,s,t).\ \)This
completes the proof.

 \end{proof}

To finish this chapter, we make the remark that for\(\ f_{t}(z)\ \)we
defined in this section we have

\[\log\left( \frac{f_{t}(z)}{e^{t}z} \right) = \sum_{k = 1}^{\infty}{c_{k}(t)}z^{k}\ \]

Recall that the logarithmic coefficients are defined as in chapter 2.3.
More precisely, if we define (See \hyperlink{D1}{[D1]} p.151 Ch 5.4)

\[\log\frac{f(z)}{z} = 2\sum_{k = 1}^{\infty}{\gamma_{k}(t)}z^{k}\]

Then \(\ c_{k}(t) = 2\gamma_{k}(t).\ \)Moreover, we may find these
coefficients by Cauchy's formula:

\[c_{k}(t) = \frac{1}{2\pi i}\int_{\gamma(0,r)}^{}{\frac{\log\left( \frac{f_{t}(z)}{e^{t}z} \right)}{z^{k + 1}}dz}\]

Although this is only a very small part of the consequences of Loewner's
theory, we are now ready to study Weinstein's paper.

\newpage

\hypertarget{proof-of-bieberbachs-conjecture-1}{%
\section{4 Proof of Bieberbach's
conjecture}\label{proof-of-bieberbachs-conjecture-1}}

\subsection{4.1 Legendre polynomials}

The notion of Legendre polynomials\(\ P_{n}(x)\) can be defined via the
generating function

\[g(x,t) = \frac{1}{\sqrt{1 - 2xt + t^{2}}} = \sum_{n = 0}^{\infty}{P_{n}(x)t^{n}}\]

It has many applications in the field of quantum mechanics and
electromagnetic theory etc. The generating function\(\ g(x,t)\ \)can be
expanded to find an explicit form of\(\ P_{n}(x).\)

Weinstein's proof relied on an important theorem concerning Legendre
polynomials, so we will spend this chapter introducing them (briefly).
Readers who are less interested in this particular field may take
Theorem 4.1.3 as given and move to the next section.

To start with, we prove the following preliminary theorems:

\uline{\textbf{Theorem 4.1.1} (Rodrigues Formula)} 

\fbox{\parbox{1.0\textwidth}{\ 

The \emph{n-th
Legendre Polynomial} is given by the following formula

\[P_{n}(x) = \frac{1}{2^{n}n!}\frac{d^{n}}{dx^{n}}\left\lbrack \left( x^{2} - 1 \right)^{n} \right\rbrack\]
}}

\begin{proof}

First notice that by expanding the generating function in series form

\[g(x,t) = \frac{1}{\sqrt{1 - 2xt + t^{2}}} = 1 + xt + \frac{1}{2}\left( 3x^{2} - 1 \right)t^{2} + \ldots\]

we have that

\[P_{0}(x) = 1,\ P_{1}(x) = x,\ P_{2}(x) = \frac{1}{2}\left( 3x^{2} - 1 \right),\ldots\]

and more generally

\[P_{n}(x) = \sum_{s = 0}^{\frac{n}{2}}{( - 1)^{s}\frac{(2n - 2s)!}{2^{n}s!(n - s)!(n - 2s)!}x^{n - 2s}}\]

for even\(\ n\ \)and

\[P_{n}(x) = \sum_{r = 0}^{\frac{n - 1}{2}}{( - 1)^{r}\frac{(2n - 2r)!}{2^{n}r!(n - r)!(n - 2r)!}x^{n - 2r}}\]

for odd\(\ n.\ \)

Now

\[\frac{d^{n}}{dx^{n}}\left( x^{2} - 1 \right)^{n} = \sum_{s = 0}^{\frac{n}{2}/\frac{n - 1}{2}}{( - 1)^{s}\frac{n!(2n - 2s)!}{s!(n - s)!(n - 2s)!}z^{n - 2s}}\]

depending on odd/even\(\ n,\ \)by Leibniz's formula. It is now easy to
see that

\[P_{n}(x) = \frac{1}{2^{n}n!}\frac{d^{n}}{dx^{n}}\left( x^{2} - 1 \right)^{n}\]

 \end{proof}

The reason why the Legendre Polynomials forms a complete (or total)
orthonormal system is because it is a result of the Gram-Schmidt process
from the spanning set\(\ \left\{ 1,x,x^{2}\ldots \right\}\ \)under the
inner product

\[\left\langle x,y \right\rangle = \int_{- 1}^{1}{x(t)y(t)dt}\ \ \]

For a detailed proof of why this guarantees completeness (which is
formally functional analysis), readers are encouraged to consult
\hyperlink{K2}{[K2]}, Theorem 3.2-3 p.139-140 and Definition 3.7-1 p.176-177. 

In view of Cauchy's formula

\[f^{n}(z) = \frac{n!}{2\pi i}\int_{\gamma}^{}{\frac{f(\xi)}{(\xi - z)^{n + 1}}d\xi}\]

It holds that

\[P_{n}(z) = \frac{1}{2^{n}n!}\frac{d^{n}}{dz^{n}}\left( z^{2} - 1 \right)^{n} = \frac{1}{2^{n}n!}\frac{n!}{2\pi i}\int_{\gamma}^{}{\frac{\left( \xi^{2} - 1 \right)^{n}}{(\xi - z)^{n + 1}}d\xi} = \frac{1}{2\pi i}\int_{\gamma}^{}{\frac{\left( \xi^{2} - 1 \right)^{n}}{{2^{n}(\xi - z)}^{n + 1}}d\xi}\]

which is referred to as \emph{Schlafli's integral formula.}

\newpage

The notion of Legendre Polynomials arises when the differential equation

\[\left( 1 - z^{2} \right)y^{''} - 2zy^{'} + \left( n^{2} + 1 \right)\lambda = 0\]

is to be solved. From Schlafli's integral formula, we see that indeed

\[\left( 1 - z^{2} \right)\left( P_{n}(z) \right)^{''} - 2z\left( P_{n}(z) \right)^{'} + \left( n^{2} + 1 \right)\lambda = \frac{(n + 1)}{2\pi i \cdot 2^{n}}\int_{\gamma}^{}{\frac{d}{d\xi}\left( \frac{\left( \xi^{2} - 1 \right)^{n + 1}}{(\xi - z)^{n + 2}} \right)}d\xi = 0\]

by the fundamental theorem of calculus and the fact that the integrand
resumes its original value after describing the closed
curve\(\ \gamma.\ \)

So the Legendre polynomials solves the above differential equation. The
same result can also be verified through substituting Rodrigue's
formula.

The key theorem that Weinstein utilized in his proof is the following:

\uline{Theorem 4.1.3 (Addition Theorem for Legendre Polynomials)}

\fbox{\parbox{1.0\textwidth}{\ 

\emph{It holds that}

\[P_{n}\left( \cos{\theta_{1}\cos{\theta_{2} + \sin{\theta_{1}\sin{\theta_{2}\cos\phi}}}} \right)\]

\emph{equals to}

\[P_{n}\left( \cos\theta_{1} \right)P_{n}\left( \cos\theta_{2} \right) + 2\sum_{k = 1}^{n}{( - 1)^{k}P_{n}^{- k}\left( \cos\theta_{1} \right)P_{n}^{k}\left( \cos\theta_{2} \right)}\cos(k\phi)\]

}}

There are various ways to prove this theorem, but none of which is very
simple. We give a comparatively simple proof to this result, based on a well-known theorem. Readers interested are
encouraged to consult \hyperlink{W2}{[W2]} p.342 Ch15.61-62 in a more detailed and rigorous manner. The above text is
also good for a further read on the theory of Legendre Polynomials.

\begin{proof}

We start by defining the \emph{Spherical Harmonics},
\(Y_{l}^{m}(\theta,\ \phi),\ \)which are closely related to Legendre
polynomials:

\[Y_{l}^{m}(\theta,\ \phi) = ( - 1)^{m}\sqrt{\frac{(2l + 1)}{4\pi}\frac{(l - m)!}{(l + m)!}}P_{l}^{m}\left( \cos\theta \right)e^{im\phi}\]

for\(\ l = 0,1,2,3\ldots\ \)and\(\ m = - l, - l + 1\ldots,l\)

Now note that for the two vectors that lies in an unit sphere

\[\mathbf{r} = \left( \cos{\phi_{1}\sin{\theta_{1},\ \sin{\phi_{1}\sin{\theta_{1},\cos\theta_{1}}}}} \right)\]

\[\mathbf{r}^{\mathbf{'}} = \left( \cos{\phi_{2}\sin{\theta_{2},\ \sin{\phi_{2}\sin{\theta_{2},\cos\theta_{2}}}}} \right)\]

it holds that

\[\mathbf{r \cdot}\mathbf{r}^{\mathbf{'}}\mathbf{=}\sin\theta_{1}\sin{\theta_{2}\left( \cos\left( \phi_{1} - \phi_{2} \right) \right) + \cos{\theta_{1}\cos\theta_{2}}}\]

\[= \sin\theta_{1}\sin{\theta_{2}\left( \cos(\phi) \right) + \cos{\theta_{1}\cos{\theta_{2}\ \ }}}\]

where\(\ \phi = \phi_{1} - \phi_{2}\). We also know that by properties
of dot products this can be view as the cosine of an angle between two
points, and name it\({\ cos}{\delta.\ }\)

From this we can derive the result

\[P_{n}\left( \cos\theta_{1} \right)P_{n}\left( \cos\theta_{2} \right) + 2\sum_{k = 1}^{n}{( - 1)^{k}P_{n}^{- k}\left( \cos\theta_{1} \right)P_{n}^{k}\left( \cos\theta_{2} \right)}\cos(k\phi)\]

\[= \frac{4\pi}{2n + 1}\sum_{m = - n}^{n}( - 1)^{m}Y_{l}^{m}\left( \theta_{1},\ \phi_{1} \right)Y_{l}^{- m}\left( \theta_{2},\ \phi_{2} \right)\]

Indeed,

\[P_{n}\left( \cos\theta_{1} \right)P_{n}\left( \cos\theta_{2} \right) + 2\sum_{k = 1}^{n}{( - 1)^{k}P_{n}^{- k}\left( \cos\theta_{1} \right)P_{n}^{k}\left( \cos\theta_{2} \right)}\cos(k\phi)\]

\[= \sum_{m = - n}^{n}( - 1)^{m}\frac{(l - m)!}{(l + m)!}P_{l}^{m}\left( \cos\theta_{1} \right)P_{l}^{m}\left( \cos\theta_{2} \right)\cos\left( k\left( \phi_{1} - \phi_{2} \right) \right)\]

\[= \sum_{m = - n}^{n}\frac{4\pi}{2l - 1}{( - 1)^{m}\left( \sqrt{\frac{(2l + 1)}{4\pi}\frac{(l - m)!}{(l + m)!}} \right)}^{2}P_{l}^{m}\left( \cos\theta_{1} \right)P_{l}^{m}\left( \cos\theta_{2} \right)\cos\left( k\left( \phi_{1} - \phi_{2} \right) \right)\]

\[= \frac{4\pi}{2n + 1}\sum_{m = - n}^{n}( - 1)^{m}Y_{l}^{m}\left( \theta_{1},\ \phi_{1} \right)Y_{l}^{- m}\left( \theta_{2},\ \phi_{2} \right)\]

by definition of Spherical Harmonics. To see that the first equation is
true, consider the fact that

\[P_{l}^{- m} = ( - 1)^{m}\frac{(l - m)!}{(l + m)!}P_{l}^{m}\left( \cos\theta_{1} \right)\]

which follows from

\[\frac{d^{l - m}}{dx^{l - m}}\left( x^{2} - 1 \right)^{l} = ( - 1)^{m}\frac{(l - m)!}{(l + m)!}\frac{d^{l + m}}{dx^{l + m}}\left( x^{2} - 1 \right)^{l}\]

More precisely, by Leibniz's formula for differentiation, we have that
(from Rodrigue's formula)

\[P_{n}^{- m} = \left\{ \frac{1}{2^{n}n!}\frac{d^{n}}{dx^{n}}\left\lbrack \left( x^{2} - 1 \right)^{n} \right\rbrack \right\}^{- m}\]

\[= \frac{( - 1)^{m}}{2^{n}n!}\left( 1 - x^{2} \right)^{\frac{m}{2}}\frac{d^{n + m}}{dx^{n + m}}\left\lbrack \left( x^{2} - 1 \right)^{n} \right\rbrack\]

\[= (x - 1)^{- \frac{m}{2}}( - 1)^{- \frac{m}{2}}(x + 1)^{- \frac{m}{2}}\sum_{j = 0}^{n - m}\begin{pmatrix}
n - m \\
j \\
\end{pmatrix}\frac{n!}{(n - j)!}(x - 1)^{n - j}\frac{n!}{(m + j)!}(x - 1)^{m + j}\]

Similarly,

\[P_{n}^{m} = (x - 1)^{\frac{m}{2}}( - 1)^{\frac{m}{2}}(x + 1)^{\frac{m}{2}}\sum_{j = m}^{n}\begin{pmatrix}
n + m \\
j \\
\end{pmatrix}\frac{n!}{(n - j)!}(x - 1)^{n - j}\frac{n!}{(j - m)!}(x - 1)^{j - m}\]

\[= ( - 1)^{\frac{m}{2}}\sum_{j = m}^{n}\frac{(n + m)!n!n!(x - 1)^{n - j + \frac{m}{2}}(x + 1)^{j - \frac{m}{2}}}{j!(n + m - j)!(n - j)!(j - m)!}\]

We use the change of index from\(\ j\ \)to\(\ k + m,\ \)such that

\[P_{n}^{m} = ( - 1)^{\frac{m}{2}}\sum_{k - 0}^{n - m}\frac{(n + m)!n!n!(x - 1)^{n - k - \frac{m}{2}}(x + 1)^{k + \frac{m}{2}}}{(k + m)!(n - k)!(n - k - m)!k!}\]

comparing this with\(\ P_{n}^{- m}\ \)yields the result.

To see that the last equation is true, notice that

\[\cos{\left( k\left( \phi_{1} - \phi_{2} \right) \right) = \frac{1}{2}\left( e^{ik\left( \phi_{1} - \phi_{2} \right)} + e^{- ik\left( \phi_{1} - \phi_{2} \right)} \right)}\]

and expand in accord with the definition of Spherical Harmonics.

Observe now that if we proof

\[P_{n}\left( \cos{\theta_{1}\cos{\theta_{2} + \sin{\theta_{1}\sin{\theta_{2}\cos\phi}}}} \right)\]

\[= P_{n}\left( \mathbf{r \cdot}\mathbf{r}^{\mathbf{'}} \right) = \frac{4\pi}{2n + 1}\sum_{m = - n}^{n}( - 1)^{m}Y_{l}^{m}\left( \theta_{1},\ \phi_{1} \right)Y_{l}^{- m}\left( \theta_{2},\ \phi_{2} \right)\]

it follows from the above result that

\[P_{n}\left( \cos{\theta_{1}\cos{\theta_{2} + \sin{\theta_{1}\sin{\theta_{2}\cos\phi}}}} \right) \]

\[= P_{n}\left( \cos\theta_{1} \right)P_{n}\left( \cos\theta_{2} \right) + 2\sum_{k = 1}^{n}{( - 1)^{k}P_{n}^{- k}\left( \cos\theta_{1} \right)P_{n}^{k}\left( \cos\theta_{2} \right)}\cos(k\phi)\]

However this is a well-known result: it is known as the
\emph{addition theorem for the spherical harmonics}. A proof can be
found in, for example, {[}A1{]} p.797.

\end{proof} 

\subsection{4.2 Remarks on proven results}

So far, we have proven enough results to understand Weinstein's proof of
Bieberbach's conjecture. This includes:

1. (From Chapter 3) There exists a family\(\ f_{t}\ \)of functions,
\(f_{t}\mathbb{:D \rightarrow C\ }\)for\(\ t \geq 0,\ \)analytic and
univalent in\(\mathbb{\ D\ }\)such that:

\[1.\ f_{0} = f;\]

\[2.\ f_{t}(z) = e^{t}z + \sum_{k = 2}^{\infty}{a_{k}(t)}z^{k}, a_{k}(t)\mathbb{\in C,\ }z \in \mathbb{D,\ }t \geq 0\]

\[3.\log\frac{f_{t}(z)}{e^{t}z} = \sum_{k = 1}^{\infty}{c_{k}(t)z^{k}}\]

\[4.\ Re\left\{ \frac{\frac{\partial f_{t}(z)}{\partial t}}{z\frac{\partial f_{t}(z)}{\partial z}} \right\} > 0 for\ all\ z \in \mathbb{D}\]

these are the properties derived from Lowner's theory: it gives the
dynamic information of the Loewner differential equation. These will be
the keystone of Weinstein's proof.

2. (From Chapter 2 the Milin conjecture)

\[\sum_{k = 1}^{n}\left( \frac{4}{k} - k\left| c_{k}(0) \right|^{2} \right)(n - k + 1) > 0,\ \ n = 1,2,3,\ldots\]

implies the Bieberbach
conjecture\(\ \left| a_{n} \right| < n\ \)for\(\ n = 2,3,\ldots\ \)

In fact, there has been thus far no proof of Bieberbach's conjecture
that does not make use of Milin's conjecture. That is, all existing
proves are proves of the Milin conjecture, which implies Bieberbach's
conjecture.

3. Theorem 4.1.3 (Addition Theorem for Legendre Polynomials)

It holds that

\[P_{n}\left( \cos{\theta_{1}\cos{\theta_{2} + \sin{\theta_{1}\sin{\theta_{2}\cos\phi}}}} \right)\]

is equivalent to

\[P_{n}\left( \cos\theta_{1} \right)P_{n}\left( \cos\theta_{2} \right) + 2\sum_{k = 1}^{n}{( - 1)^{k}P_{n}^{- k}\left( \cos\theta_{1} \right)P_{n}^{k}\left( \cos\theta_{2} \right)}\cos(k\phi)\]

We notice further that if we take\(\ \theta_{1} = \theta_{2},\ \)we have
the following lemma:

\uline{Lemma 4.2.1}

\[P_{n}\left( \cos^{2}{\theta + \sin^{2}\theta\cos\phi} \right) = \left( P_{n}\left( \cos\theta \right) \right)^{2} + 2\sum_{k = 1}^{n}{\frac{(n - k)!}{(n + k)!}\ \left( P_{n}^{k}\left( \cos\theta \right) \right)^{2}}\cos(k\phi)\]

Here we used the identity

\[P_{l}^{- m} = ( - 1)^{m}\frac{(l - m)!}{(l + m)!}P_{l}^{m}\left( \cos\theta_{1} \right)\]

With all above pieces collected, we may now finally put them together in
a remarkably short proof for Bieberbach's conjecture, which we now examine in details. 

\newpage

\subsection{4.3 Weinstein's Proof of the conjecture}

In this Chapter, we will proof Bieberbach's conjecture following
Weinstein's footpath. The readers are encouraged to read this chapter
along with {[}W1{]}, but this is not compulsory as we will quote
Weinstein's entire proof word for word.

\uline {Lemma 1} . There exists a family\(\ f_{t}\ \)of functions,
\(f_{t}\mathbb{:D \rightarrow C\ }\)for\(\ t \geq 0,\ \)analytic and
univalent in\(\mathbb{\ D\ }\)such that:

\[1.\ f_{0} = f,\]

\[2.\ f_{t}(z) = e^{t}z + \sum_{k = 2}^{\infty}{a_{k}(t)}z^{k}, \ \ \ a_{k}(t)\mathbb{\in C,\ }z \in \mathbb{D,\ }t \geq 0\]

\[3.\log\frac{f_{t}(z)}{e^{t}z} = \sum_{k = 1}^{\infty}{c_{k}(t)z^{k}},\ \ \ c_{k}(\infty) = 2/k\]

\[4.\ Re\left\{ \frac{\frac{\partial f_{t}(z)}{\partial t}}{z\frac{\partial f_{t}(z)}{\partial z}} \right\} > 0 \ \ \ for\ all\ z \in \mathbb{D}\]

\uline {Lemma 2.} (Milin conjecture)

\[\sum_{k = 1}^{n}\left( \frac{4}{k} - k\left| c_{k}(0) \right|^{2} \right)(n - k + 1) > 0,\ \ n = 1,2,3,\ldots\]

implies the Bieberbach
conjecture\(\ \left| a_{n} \right| < n\ \)for\(\ n = 2,3,\ldots\ \)

All above lemmas here has been verified and proven rigorously, and were
made as remarks in Chapter 4.2.

We now begin the proof of Bieberbach's theorem:

We aim to show that

\[\sum_{n = 1}^{\infty}\left( \sum_{k = 1}^{n}\left( \frac{4}{k} - k\left| c_{k}(0) \right|^{2} \right)(n - k + 1) \right)z^{n + 1} = \sum_{n = 1}^{\infty}{\int_{0}^{\infty}{g_{n}(t)dt\ z^{n + 1}}}\]

where\(\ g_{n}(t) \geq 0\ \)for\(\ t \geq 0,\ n = 1,2,\ldots\)
\newpage
Note that this implies

\[g_{n}(t) \geq 0 \Longrightarrow \int_{0}^{\infty}{g_{n}(t)dt\ } \geq 0 \Longrightarrow \sum_{k = 1}^{n}\left( \frac{4}{k} - k\left| c_{k}(0) \right|^{2} \right)(n - k + 1) \geq 0\]

and the scenario where equality holds will be discussed separately.

\begin{proof}

Fix\(\ z \in \mathbb{D;\ }\)define\(\ w = w_{t}(z)\ \)by

\[\frac{z}{(1 - z)^{2}} = \frac{e^{t}w}{(1 - w)^{2}}\]

with \(t \geq 0,\ \)and let\(\ z_{1} = re^{i\theta}\).

Since the Koebe function \(z/(1 - z)^{2}\) is in the class\(\ S,\ \)by
Loewner's theory there exists a parametric function\(\ w_{t}(z)\ s.t.\)
\(k(z) = \frac{z}{(1 - z)^{2}} = e^{t}k_{t}(z).\ \)This is why it is
guaranteed that the function\(\ w_{t}(z)\ \)as defined by Weinstein
indeed exists.

In particular, it holds
that\({\ w}_{t}(z) = f^{- 1}\left( e^{- t}f(z) \right),\ \)where\(\ f\ \)is
the Koebe function. Indeed,

\[\frac{e^{t}w}{(1 - w)^{2}} = \frac{z}{(1 - z)^{2}} \Longrightarrow \frac{w}{(1 - w)^{2}} = \frac{e^{- t}z}{(1 - z)^{2}} \Longrightarrow f\left( w_{t}(z) \right) = e^{- t}f(z)\]

and the result follows from taking\(\ f^{- 1}\ \)on both sides.

We are also introducing the\(\ z_{1}\ \)notion so that it differs
from\(\ z,\ \)because we will be integrating over a
variable\(\ z_{1}\ \)instead of\(\ z\ \)later in the proof.

It follows from the above point that

\[w_{t}(z) = f^{- 1}\left( e^{- t}f(z) \right) \Longrightarrow w_{0}(z) = f^{- 1}\left( e^{0}f(z) \right) = z\]

and

\[w_{t}(z) = f^{- 1}\left( e^{- t}f(z) \right) \Longrightarrow w_{t}(z) \rightarrow f^{- 1}\left( e^{- \infty}f(z) \right) = 0\ \]

as\(\ t \rightarrow \infty.\ \)For sake of convenience, we denote this
result as\(\ w_{\infty}(z) = 0.\)

We have (line 1)

\[\sum_{n = 1}^{\infty}\left( \sum_{k = 1}^{n}\left( \frac{4}{k} - k\left| c_{k}(0) \right|^{2} \right)(n - k + 1) \right)z^{n + 1}\]

(line 2)

\[= \frac{z}{(1 - z)^{2}}\sum_{k = 1}^{\infty}\left( \frac{4}{k} - k\left| c_{k}(0) \right|^{2} \right)z^{k}\]

This is by first changing the order of summation

\[\sum_{n = 1}^{\infty}\left( \sum_{k = 1}^{n}\left( \frac{4}{k} - k\left| c_{k}(0) \right|^{2} \right)(n - k + 1) \right)z^{n + 1} = \sum_{k = 1}^{\infty}{\sum_{n = k}^{\infty}\left( \frac{4}{k} - k\left| c_{k}(0) \right|^{2} \right)(n - k + 1)}z^{n + 1}\]

then replacing the index\(\ n - k + 1\ \)to\(\ m\ \)giving

\[\sum_{k = 1}^{\infty}{\sum_{n = k}^{\infty}\left( \frac{4}{k} - k\left| c_{k}(0) \right|^{2} \right)(n - k + 1)}z^{n + 1} = \sum_{k = 1}^{\infty}{\sum_{m = 1}^{\infty}\left( \frac{4}{k} - k\left| c_{k}(0) \right|^{2} \right)m}z^{m + k}\]

\[= \sum_{k = 1}^{\infty}{\sum_{m = 1}^{\infty}\left( \frac{4}{k} - k\left| c_{k}(0) \right|^{2} \right)m}z^{m}z^{k}\ \ \ \ \ \]

\[= \sum_{k = 1}^{\infty}\left( \sum_{m = 1}^{\infty}{mz^{m}} \right)\left( \frac{4}{k} - k\left| c_{k}(0) \right|^{2} \right)z^{k}\]

\[= \frac{z}{(1 - z)^{2}}\sum_{k = 1}^{\infty}{\left( \frac{4}{k} - k\left| c_{k}(0) \right|^{2} \right)z^{k}}\]

(line 3)

\[= \int_{0}^{\infty}\frac{- z}{(1 - z)^{2}}\frac{d}{dt}\left( \sum_{k = 1}^{\infty}\left( \frac{4}{k} - k\left| c_{k}(t) \right|^{2} \right)w^{k} \right)dt\]

This is by the fundamental theorem of Calculus and the fact
that\(\ w_{0}(z) = z,\ w_{\infty}(z) = 0\). Indeed, notice that
(remembering\(\ w\ \)has parameter\(\ t\))

\[\ \ (line\ 3) = \frac{z}{(1 - z)^{2}}\int_{0}^{\infty}{- \frac{d}{dt}\left( \sum_{k = 1}^{\infty}\left( \frac{4}{k} - k\left| c_{k}(t) \right|^{2} \right){w_{t}(z)}^{k} \right)}dt\]

\[ = \frac{z}{(1 - z)^{2}}\int_{\infty}^{0}{\frac{d}{dt}\left( \sum_{k = 1}^{\infty}\left( \frac{4}{k} - k\left| c_{k}(t) \right|^{2} \right){w_{t}(z)}^{k} \right)}dt\]

\[= \frac{z}{(1 - z)^{2}}\left\lbrack \sum_{k = 1}^{\infty}\left( \frac{4}{k} - k\left| c_{k}(t) \right|^{2} \right){w_{t}(z)}^{k} \right\rbrack_{\infty}^{0}\]

\[= \frac{z}{(1 - z)^{2}}\left\{ \sum_{k = 1}^{\infty}\left( \frac{4}{k} - k\left| c_{k}(0) \right|^{2} \right)w_{0}(z)^{k} - \sum_{k = 1}^{\infty}\left( \frac{4}{k} - k\left| c_{k}(\infty) \right|^{2} \right)w_{\infty}(z)^{k} \right\}\]

\[= \frac{z}{(1 - z)^{2}}\left\{ \sum_{k = 1}^{\infty}\left( \frac{4}{k} - k\left| c_{k}(0) \right|^{2} \right)z^{k} - \sum_{k = 1}^{\infty}\left( \frac{4}{k} - k\left| c_{k}(\infty) \right|^{2} \right) \cdot 0 \right\}\]

\[= \frac{z}{(1 - z)^{2}}\left\{ \sum_{k = 1}^{\infty}\left( \frac{4}{k} - k\left| c_{k}(0) \right|^{2} \right)z^{k} \right\} = line\ 2\]

where we used the results\(\ w_{0} = z,\ w_{\infty} = 0.\)

A word of notion should be added here, to prove that indeed it is
justified to say that \(\ c_{k}(\infty)\ \)is bounded and hence \(w_{\infty} = 0\ \)can be
used here to show that

\[\sum_{k = 1}^{\infty}\left( \frac{4}{k} - k\left| c_{k}(\infty) \right|^{2} \right) \cdot w_{\infty}(z)^{k} \rightarrow 0\]

To show this, notice that by definition of\(\ c_{k}(t),\ \)we may write
with Cauchy's formula

\[c_{k}(t) = \frac{1}{2\pi i}\int_{\gamma(0,r)}^{}\frac{\log\left( \frac{f_{t}(z)}{e^{t}z} \right)}{z^{k + 1}}dz\]

by the M-L inequality

\[\Longrightarrow \left| c_{k}(t) \right| \leq \frac{\max{\log\left| \frac{f_{t}(z)}{e^{t}z} \right|}}{r^{k}}\]

but by the Growth theorem (See \hyperlink{Growth}{Theorem 1.4.3})

\[\left| \frac{f_{t}(z)}{e^{t}} \right| \leq \frac{|z|}{\left( 1 - |z| \right)^{2}} \Longrightarrow \left| \frac{f_{t}(z)}{e^{t}z} \right| \leq \frac{1}{\left( 1 - |z| \right)^{2}}\]

this is bounded by a constant hence
\(\left| c_{k}(t) \right|\ \)definitely grows ``slower'' than
\(w_{t}\ \)when\(\ t \rightarrow \infty.\ \)

In Weinstein's paper it was claimed (in Lemma 1)
that\(\ c_{k}(\infty) = 2/k\). This is the case only
when\(\ f_{t}\ \)has a limit function of the Koebe function and its
rotations and not in general.

(line 4)

\[= \int_{0}^{\infty}\frac{e^{t}w}{(1 - w)^{2}}\frac{1 + w}{1 - w}\left( \sum_{k = 1}^{\infty}{k\left( c_{k}(t)c_{k}(t) \right)'}w^{k} + \sum_{k = 1}^{\infty}\left( 4 - k\left| c_{k}(t) \right|^{2} \right)w^{k}\frac{1 - w}{1 + w} \right)dt\]

To see this, note that by definition of\(\ w_{t}(z),\ \)it holds that:

\[\frac{z}{(1 - z)^{2}} = \frac{e^{t}w_{t}(z)}{\left( 1 - w_{t}(z) \right)^{2}} \Longrightarrow \frac{z}{(1 - z)^{2}}\left( 1 - w_{t}(z) \right)^{2} = e^{t}w_{t}(z)\]

Taking partial derivatives of\(\ t\ \)on both sides gives

\[\frac{\partial}{\partial t}\left\lbrack \frac{z}{(1 - z)^{2}}\left( 1 - w_{t}(z) \right)^{2} \right\rbrack = \frac{\partial}{\partial t}\left\lbrack e^{t}w_{t}(z) \right\rbrack\]

\[\Longrightarrow \frac{z}{(1 - z)^{2}}\frac{\partial}{\partial t}\left\lbrack \left( 1 - 2w + w^{2} \right) \right\rbrack = e^{t}w^{'} + e^{t}w\]

\begin{quote}
by the chain rule
\end{quote}

\[\Longrightarrow \frac{z}{(1 - z)^{2}}\left( - 2w^{'} + 2ww^{'} \right) = e^{t}w^{'} + e^{t}w\]

\[\Longrightarrow \frac{e^{t}w}{(1 - w)^{2}}\left( - 2w^{'} + 2ww^{'} \right) = e^{t}w^{'} + e^{t}w\]

\begin{quote}
by definition
\end{quote}

\[\Longrightarrow \frac{w}{(1 - w)^{2}}\left( - 2w^{'} + 2ww^{'} \right) = w^{'} + w\]

\[\Longrightarrow w^{'} = \frac{w^{2} - w}{w + 1} = w\frac{1 - w}{1 + w}\]

so

\[(line\ 3) = \int_{0}^{\infty}\frac{- z}{(1 - z)^{2}}\frac{d}{dt}\left( \sum_{k = 1}^{\infty}\left( \frac{4}{k} - k\left| c_{k}(t) \right|^{2} \right)w^{k} \right)dt\]

\[= \int_{0}^{\infty}{- \frac{e^{t}w}{(1 - w)^{2}}}\frac{d}{dt}\left( \sum_{k = 1}^{\infty}\left( \frac{4}{k} - k\left| c_{k}(t) \right|^{2} \right)w^{k} \right)dt\]

by definition of\(\ w\ \)

\[= \int_{0}^{\infty}{- \frac{e^{t}w}{(1 - w)^{2}}}\left( \sum_{k = 1}^{\infty}\left\{ \left( \frac{4}{k} - k\left| c_{k}(t) \right|^{2} \right)^{'}w^{k} + \left( \frac{4}{k} - k\left| c_{k}(t) \right|^{2} \right)\left( w^{k} \right)^{'} \right\} \right)dt\]

by the chain rule

\[= \int_{0}^{\infty}{- \frac{e^{t}w}{(1 - w)^{2}}}\left( \sum_{k = 1}^{\infty}\left\{ {- \left( kc_{k}(t)\overline{c_{k}(t)} \right)}^{'}w^{k} + \left( \frac{4}{k} - k\left| c_{k}(t) \right|^{2} \right)kw^{k}\frac{1 - w}{1 + w} \right\} \right)dt\]

by \(|z|^{2} = z \cdot \overline{z}\) and \(w'\) in terms
of\(\ w\ \)as proved above:

\[\ w^{'} = w\frac{1 - w}{1 + w} \Longrightarrow \left( w^{k} \right)^{'} = kw^{k - 1}w^{'}\]

\[\ \  = kw^{k - 1}w\frac{1 - w}{1 + w}\]

\[\ \  = kw^{k}\frac{1 - w}{1 + w}\]

Now

\[\frac{e^{t}w}{(1 - w)^{2}} = \frac{e^{t}w}{(1 - w)(1 + w)}\frac{1 + w}{1 - w} = \frac{e^{t}w}{1 - w^{2}}\frac{1 + w}{1 - w}\]

Hence after rearranging (including moving the minus sign out) we have

\[(line\ 3) = \int_{0}^{\infty}{\frac{e^{t}w}{1 - w^{2}}\frac{1 + w}{1 - w}}\left( \sum_{k = 1}^{\infty}\left\{ \left( kc_{k}(t)\overline{c_{k}(t)} \right)^{'}w^{k} + \left( 4 - k^{2}\left| c_{k}(t) \right|^{2} \right)w^{k}\frac{1 - w}{1 + w} \right\} \right)dt\]

which is line 4. Before moving to the next line, notice that since

\[\log\frac{f_{t}(z)}{e^{t}z} = \sum_{k = 1}^{\infty}{c_{k}(t)}z^{k}\]

taking partial derivative over\(\ t\ \)on both sides gives

\[\frac{e^{t}z}{f_{t}(z)}\left( \frac{e^{t}z\frac{\partial f}{\partial t} - f_{t}(z)e^{t}z}{\left( e^{t}z \right)^{2}} \right) = \sum_{k = 1}^{\infty}\frac{\partial c_{k}(t)}{\partial t}z^{k}\]

\[\  \Longrightarrow \frac{\frac{\partial f_{t}(z)}{\partial t} - f_{t}(z)}{f_{t}(z)} = \sum_{k = 1}^{\infty}\frac{\partial c_{k}(t)}{\partial t}z^{k}\]

By Cauchy's formula for Taylor coefficients,

\[\frac{\partial c_{k}(t)}{\partial t} = \frac{1}{2\pi i}\int_{|z| = r}^{}\frac{\left\lbrack \frac{\frac{\partial f_{t}(z)}{\partial t} - f_{t}(z)}{f_{t}(z)} \right\rbrack}{z^{k + 1}}dz\]

where the integral is taken over an arbitrary circle with radius\(\ r\).

Letting\(\ z = re^{i\theta}\ \)and rearranging gives

\[\frac{\partial c_{k}(t)}{\partial t} = \frac{1}{2\pi}\int_{0}^{2\pi}\frac{\frac{\frac{\partial f_{t}(z)}{\partial t}}{f_{t}(z)} - 1}{z^{k}}d\theta\]

So

\[\frac{\partial c_{k}(t)}{\partial t} = \frac{1}{2\pi}\int_{0}^{2\pi}\frac{\frac{\partial f_{t}(z)}{\partial t} - f_{t}(z)}{f_{t}(z)z^{k}}d\theta = \frac{1}{2\pi}\int_{0}^{2\pi}\frac{\frac{\partial f_{t}(z)}{\partial t}}{f_{t}(z)z^{k}}d\theta - \frac{1}{2\pi}\int_{0}^{2\pi}\frac{1}{z^{k}}d\theta\]

But as\(\ r \rightarrow 1,\ z^{- k} \rightarrow \overline{z^{k}}.\ \)
Also,

\[\frac{1}{2\pi}\int_{0}^{2\pi}\frac{1}{z^{k}}d\theta \rightarrow \frac{1}{2\pi}\int_{0}^{2\pi}e^{ik( - \theta)}d\theta = \frac{1}{2\pi}\left\lbrack - e^{ik( - \theta)} \right\rbrack_{0}^{2\pi} = 0\]

as the integral is independent of\(\ r.\ \) So we can write

\[\frac{\partial c_{k}(t)}{\partial t} = \lim_{r \rightarrow 1}\left\lbrack \frac{1}{2\pi}\int_{0}^{2\pi}\frac{\frac{\partial f_{t}(z)}{\partial t}}{f_{t}(z)}\overline{z^{k}}d\theta \right\rbrack\]

\newpage

Similarly,

\[\frac{\partial\overline{c_{k}(t)}}{\partial t} = \lim_{r \rightarrow 1}\left\lbrack \frac{1}{2\pi}\int_{0}^{2\pi}\overline{\frac{\frac{\partial f_{t}(z)}{\partial t}}{\overline{f_{t}(z)}}}z^{k}d\theta \right\rbrack\]

Now notice that

\[\frac{1}{1 - w} = \sum_{k = 0}^{\infty}w^{k} \Longrightarrow \sum_{k = 1}^{\infty}{4w^{k}} = \frac{4w}{1 - w}\]

So

\[\sum_{k = 0}^{\infty}{\left( 4 - k^{2}\left| c_{k}(t) \right|^{2} \right)w^{k}} = \sum_{k = 1}^{\infty}{4w^{k}} + \sum_{k = 1}^{\infty}{- k^{2}\left| c_{k}(t) \right|^{2}}w^{k}\]

\[= \frac{4w}{1 - w} + \sum_{k = 1}^{\infty}{- k^{2}\left| c_{k}(t) \right|^{2}}w^{k}\]

Hence line 4 to line 5-7:

\[(line\ 5-7)= \int_{0}^{\infty}{\frac{e^{t}w}{1 - w^{2}}} \left( A(w,t)\right)dt\]

Where 

\[ A(w,t)=\frac{1 + w}{1 - w}\left( 1 + \sum_{k = 1}^{\infty}\left\lbrack \lim_{r \rightarrow 1}\frac{1}{2\pi}\int_{0}^{2\pi}{\frac{\frac{\partial f_{t}\left( z_{1} \right)}{\partial t}}{f_{t}\left( z_{1} \right)}k\overline{c_{k}(t)}}\overline{z_{1}^{k}}d\theta \right\rbrack w^{k} \right) \]

\[+ \frac{1 + w}{1 - w}\left( 1 + \sum_{k = 1}^{\infty}\left\lbrack \lim_{r \rightarrow 1}\frac{1}{2\pi}\int_{0}^{2\pi}\overline{\frac{\frac{\partial f_{t}\left( z_{1} \right)}{\partial t}}{\overline{f_{t}\left( z_{1} \right)}}}kc_{k}(t)z_{1}^{k}d\theta \right\rbrack w^{k} \right)\]

\[- 2\left( \frac{1 + w}{1 - w} \right) + \frac{4w}{1 - w} + \sum_{k = 1}^{\infty}{- k^{2}\left| c_{k}(t) \right|^{2}}w^{k} \ \]

Now since

\[\frac{1 + w}{1 - w} = \frac{1 - w + 2w}{1 - w} = 1 + \frac{2w}{1 - w}\]

\[(line\ 5-7) = \int_{0}^{\infty}{\frac{e^{t}w}{1 - w^{2}}} \left( A(w,t) \right) dt\]

Where 

\[
A(w,t) = \left(1 + \frac{2w}{1 - w} \right) \left( 1 + \sum_{k=1}^{\infty} \left[ \lim_{r \rightarrow 1} \frac{1}{2\pi} \int_{0}^{2\pi} \frac{\frac{\partial f_{t}(z_{1})}{\partial t}}{f_{t}(z_{1})} k \overline{c_{k}(t)} \overline{z_{1}^{k}} \, d\theta \right] w^k \right)
\]
\[
+ \left(1 + \frac{2w}{1 - w} \right) \left( 1 + \sum_{k=1}^{\infty} \left[ \lim_{r \rightarrow 1} \frac{1}{2\pi} \int_{0}^{2\pi} \overline{\frac{\frac{\partial f_{t}(z_{1})}{\partial t}}{\overline{f_{t}(z_{1})}}} k c_{k}(t) z_{1}^{k} \, d\theta \right] w^k \right)
\]
\[
- 2 \left(1 + \frac{2w}{1 - w} \right) + \frac{4w}{1 - w} + \sum_{k=1}^{\infty} -k^{2} \left| c_{k}(t) \right|^{2} w^k
\]

Now

\[\frac{2w}{1 - w} = \sum_{j = 1}^{\infty}{2w^{j}}\]

so, denoting

\[\lim_{r \rightarrow 1}\frac{1}{2\pi}\int_{0}^{2\pi}{\frac{\frac{\partial f_{t}\left( z_{1} \right)}{\partial t}}{f_{t}\left( z_{1} \right)}k\overline{c_{k}(t)}}\overline{z_{1}^{k}}~\mathrm{d}\theta = B_{k}\]

we have that

\[\left( 1 + \frac{2w}{1 - w} \right) \left( 1 + \sum_{k = 1}^{\infty}\left\lbrack B_{k} \right\rbrack w^{k} \right)\]

\[= 1 + \sum_{k = 1}^{\infty}\left\lbrack B_{k} \right\rbrack w^{k} + 2\sum_{j = 1}^{\infty}w^{j} + 2\sum_{j = 1}^{\infty}w^{j}\sum_{k = 1}^{\infty}\left\lbrack B_{k} \right\rbrack w^{k}\]

\[= 1 + \sum_{k = 1}^{\infty}\left\lbrack B_{k} \right\rbrack w^{k} + 2\sum_{j = 1}^{\infty}w^{j} + 2\sum_{j,k = 1}^{\infty}\left\lbrack B_{k} \right\rbrack w^{k + j}\]

\[= 1 + \sum_{k = 1}^{\infty}\left\lbrack B_{k} \right\rbrack w^{k} + 2\sum_{j = 1}^{\infty}w^{j} + 2\sum_{s = 2}^{\infty}\left\lbrack B_{1} + B_{2} + \ldots + B_{s - 1} \right\rbrack w^{s}\]

\newpage

We introduce a change of index: \(\ k = s\)

\[= 1 + \sum_{s = 1}^{\infty}\left\lbrack B_{s} \right\rbrack w^{s} + 2\sum_{j = 1}^{\infty}w^{j} + 2\sum_{s = 2}^{\infty}\left\lbrack B_{1} + B_{2} + \ldots + B_{s - 1} \right\rbrack w^{s}\]

\[\ \ \ \ \ \ \ \ \ \ \ \ \ \  = 1 + B_{1}w + \sum_{s = 2}^{\infty}\left\lbrack B_{s} \right\rbrack w^{s} + 2\sum_{j = 1}^{\infty}w^{j} + 2\sum_{s = 2}^{\infty}\left\lbrack B_{1} + B_{2} + \ldots + B_{s - 1} \right\rbrack w^{s}\]

\[\  = 1 + B_{1}w + 2\sum_{j = 1}^{\infty}w^{j} + \sum_{s = 2}^{\infty}\left( 2\left\lbrack B_{1} + B_{2} + \ldots + B_{s - 1} + B_{s} \right\rbrack - B_{s} \right)w^{s}\ \ \ \ \]

by factorization. 

~\\

We introduce another change of index: \(\ j = s\)

\[= 1 + B_{1}w + 2\sum_{s = 1}^{\infty}w^{s} + \sum_{s = 2}^{\infty}\left( 2\left\lbrack B_{1} + B_{2} + \ldots + B_{s - 1} + B_{s} \right\rbrack - B_{s} \right)w^{s}\ \]

\[= 1 + B_{1}w + 2w + 2\sum_{s = 2}^{\infty}w^{s} + \sum_{s = 2}^{\infty}\left( 2\left\lbrack B_{1} + B_{2} + \ldots + B_{s - 1} + B_{s} \right\rbrack - B_{s} \right)w^{s}\]

\[\ \ \  = 1 + B_{1}w + 2w + \sum_{s = 2}^{\infty}\left( 2\left\lbrack {1 + B}_{1} + B_{2} + \ldots + B_{s - 1} + B_{s} \right\rbrack - B_{s} \right)w^{s}\]

\[= 1 + \sum_{s = 1}^{\infty}\left( 2\left\lbrack {1 + B}_{1} + B_{2} + \ldots + B_{s - 1} + B_{s} \right\rbrack - B_{s} \right)w^{s}\]

hence line 5-7 to line 8-10.

\[(line \ 8-10)= \int_{0}^{\infty}{\frac{e^{t}w}{1 - w^{2}} A(w,t) dt} \]

Where 

\[
A(w,t) = \left( 1 + \sum_{k=1}^{\infty} \left[ \lim_{r \rightarrow 1} \frac{1}{2\pi} \int_{0}^{2\pi} \frac{\frac{\partial f_{t}(z_{1})}{\partial t}}{f_{t}(z_{1})} \left( 2 \left( 1 + \ldots + k \overline{c_{k}(t) z_{1}^{k}} \right) - k \overline{c_{k}(t) z_{1}^{k}} \right) \, d\theta \right] w^{k} \right)
\]

\[+ \left( 1 + \sum_{k = 1}^{\infty}\left\lbrack \lim_{r \rightarrow 1}\frac{1}{2\pi}\int_{0}^{2\pi}\overline{\frac{\frac{\partial f_{t}\left( z_{1} \right)}{\partial t}}{\overline{f_{t}\left( z_{1} \right)}}}\left( 2\left( 1 + \ldots + kc_{k}(t)z_{1}^{k} \right) - kc_{k}(t)z_{1}^{k} \right)d\theta \right\rbrack w^{k} \right)\] 

\[- 2 + \sum_{k = 1}^{\infty}{- k^{2}\left| c_{k}(t) \right|^{2}}w^{k} \]

For the next step, notice that since

\[\log\frac{f_{t}(z)}{e^{t}z} = \sum_{k = 1}^{\infty}{c_{k}(t)}z^{k}\]

taking partial derivative over\(\ z\ \)on both sides gives

\[\ \ \ \ \ \ \ \ \frac{f_{t}^{'}(z)}{f_{t}(z)} - \frac{1}{z} = \sum_{l = 1}^{\infty}{{lc}_{l}(t)}z^{l - 1}\]

\[\Longrightarrow \frac{zf_{t}^{'}(z)}{f_{t}(z)} - 1 = \sum_{l = 1}^{\infty}{{lc}_{l}(t)}z^{l}\]

\[\Longrightarrow \frac{zf_{t}^{'}(z)}{f_{t}(z)} = 1 + \sum_{l = 1}^{\infty}{{lc}_{l}(t)}z^{l}\]

hence with some rearrangement we have from line 8-10 to 11-15 :

(line 11-15)

\[(line \ 11-15)= \int_{0}^{\infty}{\frac{e^{t}w}{1 - w^{2}}} A(w,t) dt \]

Where

\newpage

\[
A(w,t) = \sum_{k = 1}^{\infty}{- k^{2}\left| c_{k}(t) \right|^{2}}w^{k}]\]

\[ +\left( \sum_{k=1}^{\infty} \left[ \lim_{r \rightarrow 1} \frac{1}{2\pi} \int_{0}^{2\pi} \frac{\frac{\frac{\partial f_{t}(z_{1})}{\partial t}}{f_{t}(z_{1})}}{\frac{z_{1} \frac{\partial f_{t}(z_{1})}{\partial z_{1}}}{f_{t}(z_{1})}} \left\{ 1 + \sum_{l=1}^{k} {l c_{l}(t) z^{l}} \right\} \left( 2 \left( 1 + \ldots + k \overline{c_{k}(t) z_{1}^{k}} \right) - k \overline{c_{k}(t) z_{1}^{k}} \right) \, d\theta \right] w^{k} \right)
\]

\[
+ \left( \sum_{k = 1}^{\infty}\left\lbrack \lim_{r \rightarrow 1}\frac{1}{2\pi}\int_{0}^{2\pi}{\frac{\frac{\frac{\partial f_{t}\left( z_{1} \right)}{\partial t}}{\overline{f_{t}\left( z_{1} \right)}}}{\frac{z_{1}\frac{\partial f_{t}\left( z_{1} \right)}{\partial z_{1}}}{\overline{f_{t}\left( z_{1} \right)}}}\left\{ 1 + \sum_{l = 1}^{k}{{lc}_{l}(t)}z^{l} \right\}}\left( 2\left( 1 + \ldots + kc_{k}(t)z_{1}^{k} \right) - kc_{k}(t)z_{1}^{k} \right)d\theta \right\rbrack w^{k} \right)\ \]

For simplicity we denote

\[C_{0}^{k} = 1 + c_{1}(t)z_{1} + \ldots + kc_{k}(t)z_{1}^{k}\]

and

\[p\left( z_{1},t \right) = \frac{\frac{\frac{\partial f_{t}\left( z_{1} \right)}{\partial t}}{f_{t}\left( z_{1} \right)}}{\frac{z_{1}\frac{\partial f_{t}\left( z_{1} \right)}{\partial z_{1}}}{f_{t}\left( z_{1} \right)}} = \frac{\frac{\partial f_{t}\left( z_{1} \right)}{\partial t}}{z_{1}\frac{\partial f_{t}\left( z_{1} \right)}{\partial z_{1}}}\]

We expand

\[\left\{ 1 + \sum_{l = 1}^{k}{{lc}_{l}(t)}z^{l} \right\}\left( 2\left( 1 + \ldots + k\overline{c_{k}(t)z_{1}^{k}} \right) - k\overline{c_{k}(t)z_{1}^{k}} \right)\]

\[= C_{0}^{k}\left( 2\overline{C_{0}^{k}} - k\overline{c_{k}(t)z_{1}^{k}} \right)\]

\[= 2\left| C_{0}^{k} \right|^{2} - C_{0}^{k}k\overline{c_{k}(t)z_{1}^{k}}\]

\[= \frac{1}{2}\left( \left| {2C}_{0}^{k} \right|^{2} - {2C}_{0}^{k}k\overline{c_{k}(t)z_{1}^{k}} \right)\]

\[= \frac{1}{2}\left( \left| {2C}_{0}^{k} - kc_{k}(t)z_{1}^{k} \right|^{2} + 2\overline{C_{0}^{k}}kc_{k}(t)z_{1}^{k} - \left| {kc}_{k}(t)z_{1}^{k} \right|^{2} \right)\]

to see the last step note that
\(\left| {2C}_{0}^{k} - {kc}_{k}(t)z_{1}^{k} \right|^{2} = \left( {2C}_{0}^{k} - {kc}_{k}(t)z_{1}^{k} \right)\left( \overline{{2C}_{0}^{k} - kc_{k}(t)z_{1}^{k}} \right)\)

\[= \left| {2C}_{0}^{k} \right|^{2} - {2C}_{0}^{k}k\overline{c_{k}(t)z_{1}^{k}} - 2\overline{C_{0}^{k}}kc_{k}(t)z_{1}^{k} + \left| {kc}_{k}(t)z_{1}^{k} \right|^{2}\]

Similarly

\[\left\{ 1 + \sum_{l = 1}^{k}l\overline{c_{l}(t)z_{1}^{l}} \right\}\left( 2\left( 1 + \ldots + kc_{k}(t)z_{1}^{k} \right) - kc_{k}(t)z_{1}^{k} \right)\]

\[= \overline{C_{0}^{k}}\left( 2C_{0}^{k} - kc_{k}(t)z_{1}^{k} \right)\]

\[= \frac{1}{2}\left( \left| {2C}_{0}^{k} - kc_{k}(t)z_{1}^{k} \right|^{2} + 2C_{0}^{k}\overline{{kc}_{k}(t)z_{1}^{k}} - \left| {kc}_{k}(t)z_{1}^{k} \right|^{2} \right)\]

We may now re-arrange line 11-15 to the form

\[= \int_{0}^{\infty}{\frac{e^{t}w}{1 - w^{2}}} A(w,t) dt \]

Where

\[
A(w,t) = \sum_{k=1}^{\infty} -k^{2} \left| c_{k}(t) \right|^{2} w^{k}\]

\[+\left( \sum_{k=1}^{\infty} \left[ \lim_{r \rightarrow 1} 2 \, \text{Re} \left\{ \frac{1}{2\pi} \int_{0}^{2\pi} p\left( z_{1}, t \right) \frac{1}{2} \left( \left| 2C_{0}^{k} - k c_{k}(t) z_{1}^{k} \right|^{2} + 2 \overline{C_{0}^{k}} k c_{k}(t) z_{1}^{k} - \left| k c_{k}(t) z_{1}^{k} \right|^{2} \right) d\theta \right\} \right] w^{k} \right) \]

Now

\[2 \times \frac{1}{2\pi}\int_{0}^{2\pi}{p\left( z_{1},t \right)}\frac{1}{2}\left( \left| {2C}_{0}^{k} - kc_{k}(t)z_{1}^{k} \right|^{2} + 2\overline{C_{0}^{k}}kc_{k}(t)z_{1}^{k} - \left| {kc}_{k}(t)z_{1}^{k} \right|^{2} \right)d\theta \]

\[= \frac{1}{2\pi}\int_{0}^{2\pi}{p\left( z_{1},t \right)}\left| {2C}_{0}^{k} - kc_{k}(t)z_{1}^{k} \right|^{2}d\theta + 2 \times \frac{1}{2\pi}\int_{0}^{2\pi}{p\left( z_{1},t \right)}\overline{C_{0}^{k}}kc_{k}(t)z_{1}^{k}d\theta \]

\[- \frac{1}{2\pi}\int_{0}^{2\pi}{p\left( z_{1},t \right)}\left| {kc}_{k}(t)z_{1}^{k} \right|^{2}d\theta\]

But

\[\overline{C_{0}^{k}}kc_{k}(t)z_{1}^{k} = kc_{k}(t)z_{1}^{k} + kc_{k}(t)\sum_{n = 1}^{k - 1}{nr^{2n}\overline{c_{n}(t)}z_{1}^{k - n}} + k^{2}\left| c_{k}(t) \right|^{2}\left| z_{1} \right|^{2k}\]

and using the identity

\[g(a) = \frac{1}{2\pi r}\int_{0}^{2\pi}{g\left( a + re^{i\theta} \right)}d\theta\]

(where\(\ g\ \)is any analytic function on the closed disk\(\ D(a,r)\))

we have that

\[2 \times \frac{1}{2\pi}\int_{0}^{2\pi}{p\left( z_{1},t \right)}\overline{C_{0}^{k}}kc_{k}(t)z_{1}^{k}d\theta \]

\[= 2p(0,t)\left( {kc}_{k}(t) \times 0 + kc_{k}(t)\sum_{n = 1}^{k - 1}{nr^{2n}\overline{c_{n}(t)} \times 0} + k^{2}\left| c_{k}(t) \right|^{2}\left| z_{1} \right|^{2k} \right)\]

\[= p(0,t)k^{2}\left| c_{k}(t) \right|^{2}\left| z_{1} \right|^{2k}\ \ \ \ \ \ \ \ \ \ \ \ \ \ \ \ \ \]

Similarly

\[\frac{1}{2\pi}\int_{0}^{2\pi}{p\left( z_{1},t \right)}\left| {kc}_{k}(t)z_{1}^{k} \right|^{2}d\theta = p(0,t)k^{2}\left| c_{k}(t) \right|^{2}\left| z_{1} \right|^{2k}\ \ \]

Now by definition of\(\ f_{t},\ \)we have that for all\(\ t,\ \)

\[p(0,t) = \frac{e^{t}}{e^{t}} = 1\]

And hence

\[\lim_{r \rightarrow 1}\ 2\ \left\{ \frac{1}{2\pi}\int_{0}^{2\pi}{p\left( z_{1},t \right)}\frac{1}{2}\left( \left| {2C}_{0}^{k} - kc_{k}(t)z_{1}^{k} \right|^{2} + 2\overline{C_{0}^{k}}kc_{k}(t)z_{1}^{k} - \left| {kc}_{k}(t)z_{1}^{k} \right|^{2} \right)d\theta \right\}\]

\[= \frac{1}{2\pi}\int_{0}^{2\pi}{p\left( z_{1},t \right)}\left| {2C}_{0}^{k} - kc_{k}(t)z_{1}^{k} \right|^{2}d\theta + k^{2}\left| c_{k}(t) \right|^{2}\]

Substituting into line 11-15 and rearranging we have that

\newpage

(line 16-17)

\[= \int_{0}^{\infty}{\frac{e^{t}w}{1 - w^{2}}\sum_{k = 1}^{\infty}{\lim_{r \rightarrow 1}\left( \frac{1}{2\pi}\int_{0}^{2\pi}{Re\ \left\{ {p\left( z_{1},t \right)\left| {2C}_{0}^{k} - kc_{k}(t)z_{1}^{k} \right|}^{2}d\theta \right\}} \right)}w^{k}dt}\]

and if we write

\[A_{k}(t) = \lim_{r \rightarrow 1}\left( \frac{1}{2\pi}\int_{0}^{2\pi}{Re\ \left\{ p\left( z_{1},t \right) \right\}}d\theta{p\left( z_{1},t \right)\left| {2C}_{0}^{k} - kc_{k}(t)z_{1}^{k} \right|}^{2}d\theta \right) \geq 0\]

(by the fact that\(\ Re\left( p\left( z_{1},t \right) \right) > 0\) and
squared constants are also non-negative)

Then we arrive at

(line 18)

\[= \int_{0}^{\infty}{\frac{e^{t}w}{1 - w^{2}}\left( \sum_{k = 1}^{\infty}{A_{k}(t)}w^{k} \right)dt}\]

~\\

We now show that

\[\frac{e^{t}w}{1 - w^{2}} = \sum_{n = 0}^{\infty}{\Lambda_{k}^{n}(t)z^{n + 1}}\]

where\(\ \Lambda_{k}^{n}(t) \geq 0\ \)for\(\ t \geq 0.\ \) Recall from Theorem 4.1.3 (\emph{Addition theorem for Legendre Polynomials}) that

\[P_{n}\left( \cos{\theta_{1}\cos{\theta_{2} + \sin{\theta_{1}\sin{\theta_{2}\cos\phi}}}} \right)\]

is equivalent to

\[P_{n}\left( \cos\theta_{1} \right)P_{n}\left( \cos\theta_{2} \right) + 2\sum_{k = 1}^{n}{( - 1)^{k}P_{n}^{- k}\left( \cos\theta_{1} \right)P_{n}^{k}\left( \cos\theta_{2} \right)}\cos(k\phi)\]

If we let
\(\cos\theta_{1} = \cos{\theta_{2} = \left( 1 - e^{- t} \right)^{1/2}},\ \)and
write

\[\cos{\theta_{1}\cos{\theta_{2} + \sin{\theta_{1}\sin{\theta_{2}\cos\phi}}}} = 1 - e^{- t} + e^{- t}\cos\phi\]

setting

\[\cos{\delta = 1 - e^{- t} + e^{- t}\cos\phi}\]

Then it is easy to see that, from definition of\(\ w_{t},\ \)(recall
from the start of the proof)

\[\frac{z}{1 - 2\cos\delta z + z^{2}} = \frac{1}{z + \frac{1}{z} - 2\cos\delta} = \frac{1}{2 + e^{- t}\left( \frac{1}{w} + w - 2 \right) - 2\cos\delta}\]

here we used the fact that

\[\frac{z}{(1 - z)^{2}} = \frac{e^{t}w}{(1 - w)^{2}} \Longrightarrow \frac{(1 - z)^{2}}{z} = \left( \frac{(1 - w)^{2}}{e^{t}w} \right) \Longrightarrow z + \frac{1}{z} - 2 = e^{- t}\left( \frac{1 + w^{2} - 2w}{w} \right)\]

Multiplying both the denominator and numerator by\(\ e^{t}w\ \)gives

\[\ \frac{z}{1 - 2\cos\delta z + z^{2}} = \frac{e^{t}w}{2e^{t}w + \left( 1 + w^{2} - 2w \right) - 2\cos{\delta e^{t}w}}\]

\[= \frac{e^{t}w}{2e^{t}w + \left( 1 + w^{2} - 2w \right) - 2\left( 1 - e^{- t} + e^{- t}\cos\phi \right)e^{t}w}\]

\[= \frac{e^{t}w}{1 - 2w\cos{\phi + w^{2}}}\]

\[= \frac{e^{t}w}{1 - w^{2}}\frac{1 - w^{2}}{1 - 2w\cos{\phi + w^{2}}}\]

\[= \frac{e^{t}w}{1 - w^{2}} \frac {1}{2} \left(\frac{1 + e^{i\phi}w}{1 - e^{i\phi}w} + \frac{1 + e^{-i\phi}w}{1 - e^{-i\phi}w}\right)\]

\[= \frac{e^{t}w}{1 - w^{2}}\ Re\left( \frac{1 + e^{i\phi}w}{1 - e^{i\phi}w} \right)\]

\[= \frac{e^{t}w}{1 - w^{2}}\ Re\left( \frac{1 - e^{i\phi}w + 2e^{i\phi}w}{1 - e^{i\phi}w} \right)\]

\[= \frac{e^{t}w}{1 - w^{2}}\ Re\left( 1 + 2\frac{e^{i\phi}w}{1 - e^{i\phi}w} \right)\]

\[= \frac{e^{t}w}{1 - w^{2}}\left( 1 + 2\sum_{k = 1}^{\infty}{w^{k}\cos{k\phi}} \right)\]

\[= \frac{e^{t}w}{1 - w^{2}} + 2\sum_{k = 1}^{\infty}{{\sum_{n = 1}^{\infty}{\Lambda_{k}^{n}(t)z^{n + 1}}\cos}{k\phi}}\]

since we defined

\[\sum_{n = 1}^{\infty}{\Lambda_{k}^{n}(t)z^{n + 1}} = \frac{e^{t}w}{1 - w^{2}}\]

so

\[\frac{z}{1 - 2\cos\delta z + z^{2}} = \frac{e^{t}w}{1 - w^{2}} + 2\sum_{k = 1}^{\infty}{{\sum_{n = 1}^{\infty}{\Lambda_{k}^{n}(t)z^{n + 1}}\cos}{k\theta}}\]

But notice that

\[\frac{1}{\sqrt{1 - 2\cos\delta z + z^{2}}}\]

is the generating function for the Legendre Polynomials. Indeed, recall
from the definition of Legendre Polynomials from Chapter 4.1,

\[g(x,t) = \frac{1}{\sqrt{1 - 2xt + t^{2}}} = \sum_{n = 0}^{\infty}{P_{n}(x)t^{n}}\]

so

\[\frac{1}{\sqrt{1 - 2\cos\delta z + z^{2}}} = \sum_{n = 0}^{\infty}{P_{n}\left( \cos\delta \right)z^{n}}\]

and from Lemma 4.1.2 we have that
when\({\ cos}{\left( \theta_{1} \right) =}{\ cos}\left( \theta_{2} \right)\),

\[P_{n}\left( \cos\delta \right) = \left( P_{n}\left( \cos\theta \right) \right)^{2} + 2\sum_{k = 1}^{n}{\frac{(n - k)!}{(n + k)!}\ \left( P_{n}^{k}\left( \cos\theta \right) \right)^{2}}\cos(k\phi)\]

Thus

\[\ \ \ \frac{1}{\sqrt{1 - 2\cos\delta z + z^{2}}} = \sum_{n = 0}^{\infty}{P_{n}\left( \cos\delta \right)z^{n}}\]

\[= \sum_{n = 0}^{\infty}{\left( P_{n}\left( \cos\theta \right) \right)^{2}z^{n} + 2}\sum_{k = 1}^{\infty}{\left( \sum_{n = k}^{\infty}{\frac{(n - k)!}{(n + k)!}\ \left( P_{n}^{k}\left( \cos\theta \right) \right)^{2}}z^{n} \right)\cos{k\phi}}\]

and the coefficients are squared real constants, hence non-negative. It
only remains to see that

\[\frac{z}{1 - 2\cos\delta z + z^{2}} = \left( \frac{1}{\sqrt{1 - 2\cos\delta z + z^{2}}} \right)^{2}z = \frac{e^{t}w}{1 - w^{2}} + 2\sum_{k = 1}^{\infty}{{\sum_{n = 1}^{\infty}{\Lambda_{k}^{n}(t)z^{n + 1}}\cos}{k\theta}}\]

so obviously \(\Lambda_{k}^{n}(t) \geq 0\ \)for all\(\ t \geq 0.\ \)

With the above argument, we have established the fact that

\[\sum_{n = 1}^{\infty}\left( \sum_{k = 1}^{n}\left( \frac{4}{k} - k\left| c_{k}(0) \right|^{2} \right)(n - k + 1) \right)z^{n + 1} = \int_{0}^{\infty}{\frac{e^{t}w}{1 - w^{2}}\left( \sum_{k = 1}^{\infty}{A_{k}(t)}w^{k} \right)dt}\]

\[= \int_{0}^{\infty}{\left( \sum_{n = 1}^{\infty}{\Lambda_{k}^{n}(t)z^{n + 1}} \right)\left( \sum_{k = 1}^{\infty}{A_{k}(t)}w^{k} \right)dt} = \sum_{n = 1}^{\infty}{\int_{0}^{\infty}{g_{n}(t)dt\ z^{n + 1}}}\]

where\(\ g_{n}(t) \geq 0.\ \)This implies the Milin conjecture as
required, and, henceforth, proves Bieberbach's conjecture.

\end{proof}

\newpage

\section{Footnote}\label{footnote}

This paper was started in December 2022, term 1 of my
1\textsuperscript{st} year as an undergrad, when I was first introduced
to univalent function theory. Difficulties encountered during the study
of this field were largely due to the fact that most related reading
materials faces advanced readers, sketching or skipping proves that are
challenging for elementary audiences. The purpose of this paper is
fulfilled, have it interested or aided readers in making a start in the
study of geometric function theory.

I would like to thank my supervisor (Prof. A.Sobolev) again, for
introducing me to this particular area, and for encouraging me during
hard times: an understanding that values too much for me to put in
words. I would also like to thank M.P, R.Reynolds, T.Ida, E.Cook, R.Zhang, my
parents and all other personnel who supported me throughout this
project.

\newpage

\hypertarget{references}{%
\section{References}\label{references}}

\hypertarget{A1}{[}A1{]} George B. Arfken, Hans J. Weber etc, \emph{Mathematical Methods for Physicists, Seventh Edition A Comprehensive Guide, ISBN 9780123846549}

~

\hypertarget{A2}{[}A2{]} Andrew Presley, \emph{Elementary Differential Geometry,
Springer Undergraduate Mathematics Series, ISBN 978-1848828902}

~

\hypertarget{B1}{[}B1{]}Bak and Newman, \emph{Complex Analysis, Undergraduate texts in
mathematics, 3\textsuperscript{rd} edition (06 August 2010), ISBN
978-1-4419-7287-3}

~

\hypertarget{Br1}{[}Br1{]} Louis de Branges. \emph{A proof of the Bieberbach conjecture} Acta Math. 154 (1-2) 137 - 152, 1985. https://doi.org/10.1007/BF02392821

~

\hypertarget{C1}{[}C1{]} Conway, \emph{Complex analysis in the spirit of Lipman Bers,
ISBN 978-0-387-74714-9}

~

\hypertarget{D1}{[}D1{]} Duren, \emph{Univalent Functions, ISBN 9780387907956}, 1983

~

\hypertarget{P1}{[}P1{]} Pommerenke, \emph{Univalent Functions (with a chapter on
Quadratic Differentials by Gerd Jensen)} , ISBN 3-525-40133-7, 1975

~

\hypertarget{K1}{[}K1{]} Antti Kemppainen, \emph{Schramm--Loewner Evolution, ISBN
978-3-319-65327-3}

~

\hypertarget{K2}{[}K2{]} Erwin Kreyszig, \emph{Introductory functional analysis with
applications, ISBN 0-471-50731-8}

~

\hypertarget{KO1}{[}KO1{]} Koepf, W. \emph{Bieberbach’s conjecture, the de Branges and Weinstein functions and the Askey-Gasper inequality. Ramanujan J 13, 103–129 (2007). https://doi.org/10.1007/s11139-006-0244-2}

~

\hypertarget{S1}{[}S1{]} Barry Simon, \emph{Advanced Complex Analysis A Comprehensive
Course in Analysis, Part 2B, ISBN 9781470411015}

~

\hypertarget{S2}{[}S2{]} Elias M. Stein, \emph{Complex Analysis ISBN 978-0691113852}

~

\hypertarget{W1}{[}W1{]} Weinstein, \emph{The Bieberbach Conjecture}, International
Mathematics Research Notices, 1991

~

\hypertarget{W2}{[}W2{]} E. T. Whittaker, G. N. Watson, \emph{A Course of Modern
Analysis,}

~

\hypertarget{Wi1}{[}Wi1{]} Wilf, \emph{A Footnote on Two Proofs of the Bieberbach-De Branges Theorem, Bulletin of the London Mathematical Society, Volume 26, Issue 1, January 1994, Pages 61–63, https://doi.org/10.1112/blms/26.1.61}

\emph{ISBN 9781316518939}

\end{document}